\documentclass[12pt]{aptpub}
\usepackage{amsfonts,amsmath,latexsym,amssymb,mathrsfs}
\usepackage{graphicx}

\numberwithin{equation}{section}

\def\numberlikeadb{\global\def\theequation{\thesection.\arabic{equation}}}
\numberlikeadb

\newcommand{\RR}{{\bf R}}

\newcommand{\beas}{\begin{eqnarray*}}
\newcommand{\enas}{\end{eqnarray*}}
\newcommand{\eqs}{\begin{eqnarray*}}
\newcommand{\ens}{\end{eqnarray*}}

\newcommand{\bea}{\begin{eqnarray}}
\newcommand{\eqa}{\begin{eqnarray}}
\newcommand{\ena}{\end{eqnarray}}
\newcommand{\eq}{\begin{equation}}
\newcommand{\en}{\end{equation}}

\newcommand{\proofbox}{\hspace*{\fill}\mbox{$\halmos$}}
\newcommand{\halmos}{\rule{1ex}{1.4ex}}

\def\ignore#1{}

\def\half{{\textstyle{\frac12}}}
\def\shalf{{\scriptstyle{\frac12}}}
\def\stq{{\scriptstyle{\frac34}}}
\def\third{{\textstyle{\frac13}}}
\def\twothirds{{\textstyle{\frac23}}}
\def\quarter{{\textstyle{\frac14}}}
\def\eighth{{\textstyle{\frac18}}}

\def\Ref#1{(\ref{#1})}
\def\a{\alpha}
\def\b{\beta}
\def\s{\sigma}
\def\f{\phi}

\def\l{\lambda}
\def\L{\Lambda}

\def\dtv{d_{TV}}
\def\law{{\cal L}}

\def\ep{\hfill $\proofbox$ \bigskip}
\def\Def{\ :=\ }

\def\re{\RR}

\def\giv{\,|\,}

\def\Po{{\rm Po\,}}

\def\non{\nonumber}
\def\th{\theta}

\def\e{\varepsilon}
\def\m{\mu}

\def\var{{\rm Var\,}}

\def\r{\rho}

\def\t{\tau}
\def\g{\gamma}
\def\h{\eta}
\def\f{\phi}

\def\ps{\psi}

\def\lti{{\lim_{t\to\infty}}}

\def\Bl{\left(}
\def\Br{\right)}

\def\Blb{\left\{}
\def\Brb{\right\}}

\def\ui{^{(1)}}
\def\ut{^{(2)}}

\def\Bi{{\rm Bi\,}}

\def\ex{{\mathbb E}}
\def\pr{{\mathbb P}}

\def\ff{{\cal F}}
\def\d{\delta}

\def\uh{^{(3)}}

\def\Eq{\ =\ }
\def\Le{\ \le\ }

\def\sjo{\sum_{j\ge0}}

\def\adb{{}}

\def\tH{{\widetilde H}}

\def\ignore#1{}

\def\uj{^{(j)}}

\def\n{\nu}

\def\Giv{\,\Big|\,}
\def\cupdot{\cup\kern-8.2pt\cdot\kern5.5pt}

\def\tc{{\tilde c}}

\def\bone{{\bf{1}}}

\def\tU{{\widetilde U}}

\def\tW{{\widetilde W}}

\def\Def{\ :=\ }

\def\MBP{X^*}

\def\YYY{Y}

\def\tM{{\widetilde M}}
\def\KK{{\cal K}}

\def\p{\pi}
\def\iud{independently and uniformly distributed}
\def\tit{{\tilde\t}}
\def\tP{{\widetilde P}}
\def\tC{{\widetilde C}}
\def\bC{{\overline \YYY}}
\def\tQ{{\widetilde Q}}
\def\o{\omega}
\def\sjd{\sum_{j=0}^d}
\def\sji{\sum_{j\ge1}}
\def\mm{{\mathfrak m}}
\def\tZ{{\widetilde Z}}
\def\tA{{\widetilde A}}
\def\ppp{q}
\def\YYYs{\YYY^*}
\def\uip{^{(1)+}}
\def\utp{^{(2)+}}
\def\uhp{^{(3)+}}
\def\ulp{^{(l)+}}
\def\Xm{{X^-}}
\def\Xp{{X^+}}
\def\exs{\ex_{(s)}}
\def\ta{{\tilde a}}
\def\tU{{\widetilde U}}

\authornames{A. D. Barbour and G. Reinert}
\shorttitle{Gossip processes}

\begin{document}

\title{Asymptotic behaviour of gossip processes and small world networks}

\authorone[Universit\"at Z\"urich and National University of Singapore]{A. D. Barbour}
\addressone{Angewandte Mathematik, Universit\"at Z\"urich,
Winterthurertrasse 190, CH-8057 Z\"URICH;  ADB was Visiting Research Professor
at the National University of Singapore while part of this work was carried out.
Work supported in part by Australian Research Council Grants Nos DP120102728 and DP120102398}

\authortwo[University of Oxford]
{G. Reinert}
\addresstwo{Department of Statistics,
University of Oxford, 1 South Parks Road, OXFORD OX1 3TG, UK.
GDR was supported in part by EPSRC and BBSRC through OCISB.
}


\begin{abstract}
Both small world models of random networks with occasional long range connections and
gossip processes with occasional long range transmission of information have
similar characteristic behaviour.  The long range elements appreciably reduce the
effective distances, measured in space or in time, between pairs of typical points.
In this paper, we show that their common behaviour can be interpreted as a product of
the locally branching nature of the models.  In particular, it is shown that both typical
distances between points and the proportion of space that can be reached within a given
distance or time can be approximated by formulae involving the limit random variable
of the branching process.
\end{abstract}


\section{Introduction}\label{intro}
 \setcounter{equation}{0}

Moore \& Newman~\cite{MooreNewman} introduced a continuous analogue of the Watts
\& Strogatz~\cite{WattsStrogatz} ``small world'' model.  In this model, a random
number of chords are superimposed as shortcuts on a circle~$C$ of
circumference~$L$. The chords have endpoints uniformly and independently distributed
on~$C$, and the number of chords follows a Poisson distribution $\Po(L\rho/2)$ with mean~$L\rho/2$,
for some $\rho=\rho(L)$.  Distance is measured as usual along the
circumference, but chords are deemed to be of length zero, and
interest centres on finding the statistics of shortest path
distances between pairs of points.  A closely related model, the
``great circle model'', was introduced somewhat earlier by Ball {\it
et al.}\ \cite{Balletal}, in the context of epidemics; here, distance
between points translates into time taken for one infected person to
infect another. In~Barbour \& Reinert~\cite{BR}, assuming the
expected number $L\rho/2$ of shortcuts to be large, we proved a
distributional approximation for the distance between a randomly
chosen pair of points $P$ and~$P'$, and gave a bound on the order of
the error, in terms of total variation distance. We also showed that
analogous results could be proved in higher dimensions by much the
same method, when the circle is replaced by a sphere or a torus.  It
turns out that the reduction in the typical distances between pairs
of points that results from introducing shortcuts is still
substantial, but less dramatic than in one dimension.

More recently,  Chatterjee \& Durrett~\cite{ChatDurr} studied a model
for the spread of gossip that is the continuous analogue of one of a number of models
discussed in Aldous~\cite{Aldous}.  Here, information spreads locally from
an individual to his neighbours on the two-dimensional torus, and also occasionally
to other, randomly
chosen members of the community.  Thus a disc of informed individuals,
centred on an initial informant, grows steadily in the torus;
long range transmissions of information occur in a Poisson process, whose
rate is proportional to the area (number) of informed individuals, and any
such transmission contacts a randomly chosen point of the torus, initiating
a new disc of informed individuals.  The distinction between this model
and the corresponding two dimensional model in~\cite{BR} is that, in the gossip model,
the Poisson process runs at a rate proportional to the area of the currently
informed region; in~\cite{BR}, where the Poisson number of shortcuts is considered
to be fixed in advance, the Poisson process corresponds to a process of discovery of
shortcuts, and its rate is thus proportional to the length of the boundary
of the informed region.

Here, we consider
the development of such a process~$\YYY$ on a smooth closed homogeneous Riemannian manifold~$C$ of dimension~$d$,
such as a sphere or a torus, having large finite volume~$|C| =: L$ with respect to its intrinsic
metric.  We assume that, around each point~$P$
of~$C$, there is a collection of closed subsets~$\KK(P,s)$, $s\ge0$, that are balls of radius~$s$
with respect to a metric $d_C$ that makes~$C$ a geodesic space, and with (intrinsic) volumes $v_s(\KK) := |\KK(P,s)|
\sim s^d v(\KK)$ as $s \to 0$, for some $v(\KK) > 0$; $s$ is thought of as time, and $\{v(\KK)\}^{1/d}$
as a (linear) speed of propagation.
The metric $d_C$ need not be the same as the intrinsic metric; for instance,  on the torus, we could consider rectangular as well as circular
neighbourhoods.
The set~$\KK(P,s)$ denotes the set of points `locally' contacted after
time~$s$ has elapsed following an initial `long range' contact at~$P$, thought of as
`islands' in~$C$, and the complete set
of contacts~$\YYY_{P_0}(t)$ at time~$t$ is the union of these sets growing from an initial
point~$P_0$, and from all long range contacts made before~$t$;
we denote its volume by~$V_{P_0}(t)$.  The rate at
which long range contacts are made is proportional either to the area of the boundary
of~$\YYY_{P_0}(t)$ (small world processes) or to its volume (gossip processes), and we denote the
constant of proportionality by~$\r$;
long range contacts are made to independently and uniformly chosen points of~$C$.

The main result of this paper is Theorem~\ref{path-theorem}, which establishes a
pathwise law of large numbers, together with a rough error bound, for the time evolution
of the covered fraction $L^{-1}V_{P_0}(t)$ of~$C$, in the setting of a quite
general gossip process.
An analogous result is stated for small world processes in Theorem~\ref{path-theorem-2}.
Theorem~\ref{path-theorem} shows that there exists a random variable~$U$ such that,
for positive constants $a_1,a_2$ and $c<\infty$,
\[
   \pr\Bigl[\sup_x|L^{-1}V_{P_0}(\l_0^{-1}\{\log\L + x\}) - h_d(x + \log C_d + U)| > 4\L^{-a_1}\Bigr]
    \Le c\L^{-a_2}.
\]
Here, $C_d = (d+1)^{-1}d!$, the function~$h_d$ depends
only on the dimension~$d$, $\l_0$ is the initial exponential
growth rate of the process~$Y_{P_0}$, and~$\L := L\l_0^d/v(\KK)$.
The value taken by the random variable~$U$ is essentially
determined by the very early evolution of~$Y_{P_0}$, and can be
thought of as a random delay, caused by early fluctuations in the
growth of the process, before the deterministic evolution governed
by~$h_d$ sets in.  The function $h_d$ is defined through a Laplace
transform and satisfies an integral equation, \eqref{h-eqn}.
Both~$U$ and~$h_d$ have interpretations in terms of an associated
Markov branching process~$\MBP$.  In the particular case of the
torus, our result extends the limit law proved by Chatterjee \&
Durrett~\cite{ChatDurr}, by providing an estimate of the
approximation error that is uniform for all time.

Our argument is developed from that in~\cite{BR}.
The key observation is that, at least for a while, the process~$Y_{P_0}$ can
be closely approximated using a Markovian growth and branching process~$\MBP = \MBP_{P_0}$,
and that this approximation is accurate enough for the calculations that need to be made.
The state~$\MBP(t)$ of the Markov branching process at time~$t$ consists of $\KK(P_0,t)$,
together with a collection
of some number~$n(t)$ of sets of the form $\KK(P_j^*,t-\t_j^*)$, where $0\le \t_j^* \le t$
and $P_j^* \in C$ for each $1\le j\le n(t)$.  The~$P_j^*$ are independently and
uniformly chosen points of~$C$, and arise as the points~$\t_j^*$ of a Poisson process,
whose rate depends on the current state of~$\MBP$.
In `gossip' models, new contacts are made at a rate proportional to
the current volume, which, for the process~$\MBP$, is given by
\[
    V^*(t) \Def \sum_{j=0}^{n(t)} v_{t-\t^*_j}(\KK)
   \ \sim\ \sum_{j=0}^{n(t)} (t-\t^*_j)^{d} v(\KK);
\]
in `small world' models, the rate is proportional to the derivative of
the volume. The set $\YYY_{P_0}(t)$
can initially be approximated by the union
\eq\label{union-state}
   \YYYs_{P_0}(t) \Def \bigcup_{j=0}^{n(t)} \KK(P_j^*,t-\t_j^*),
\en
where $\t_0^* := 0$ and $P_0^* = P_0$. Indeed, one can take $\YYY_{P_0}(t) = \YYYs_{P_0}(t)$ until the
(random) time~$\widehat T$ at which the union in~\Ref{union-state} ceases to be
disjoint.  Thereafter,
the rate of contacts is smaller in~$\YYY_{P_0}$ than it is in~$\YYYs_{P_0}$,
and the two processes progressively separate.

In \cite{BR}, the distribution of inter-point distances
in the small world model is determined by running
two such branching processes from randomly chosen initial points $P_0'$ and~$P_0''$, each for
a length of time~$t^*$ at which the mean number of overlaps in~\Ref{union-state}
is of order~$O(1)$.  At this time, conditionally on the contact times in
the two branching processes, the number of {\it permissible\/} overlaps  between the
sets $Y_{P'_0}^*$ and~$Y_{P''_0}^*$ --- cases in which an island in one branching process
is contained within an island in the other could not have arisen in the actual
small world
process --- has an approximately Poisson distribution, and the distance between the
initial points is greater than~$2t^*$ if there are no permissible overlaps.  In this way,
and by varying the choice of~$t^*$ appropriately, the
distribution of inter-point distances can be approximated, without ever having
to go into the dependence structure that becomes important in the
process~$\YYY$ at times significantly larger than~$t^*$.  In contrast,
Chatterjee \& Durrett~\cite{ChatDurr} go beyond the branching phase in the analysis
of~$\YYY$ in their two-dimensional gossip model, and are able to prove a
conditional law of large numbers for the fraction of the torus contained in~$\YYY$,
given the outcome of the branching phase.  They also establish the asymptotics
of the first time at which $\YYY$ covers the whole torus.

The point of departure for our argument is that $\ex\{ V_{P_0}(t)/L \giv \ff_s\}$,
the conditional expectation of the covered fraction at~$t$, given the history
of~$\YYY_{P_0}$ up to time~$s$, including the initial point~$P_0$, is given by
\[
   \ex\{ V_{P_0}(t)/L \giv \ff_s\} \Eq \pr[Q \in \YYY_{P_0}(t) \giv \ff_s],
\]
where~$Q$ denotes an independently and uniformly distributed point of~$C$.
Now, for~$t$ in the relevant range (corresponding to~$2t^*$ above),
we have
\eq\label{set-intersection}
    Q \in \YYY_{P_0} \quad \mbox{if and only if} \quad \YYY_{P_0}(t/2)
      \cap \overline\YYY_Q(t/2) \neq \emptyset,
\en
where~$\overline\YYY_Q$ is an independent gossip process started from~$Q$.
The probability of the latter event can then be closely enough approximated
by computing the probability of intersection of independent branching processes
$\MBP_{P_0}$ and~$\MBP_Q$ at time~$t/2$.

The law of large numbers is proved by using an argument of much the
same flavour, since
\[
  \ex\{ [V_{P_0}(t)/L]^2 \giv \ff_s\} \Eq
    \pr[\{Q \in \YYY_{P_0}(t)\} \cap \{Q' \in \YYY_{P_0}(t)\} \giv \ff_s],
\]
for two independent and uniformly distributed points $Q,Q' \in C$.
Using~\Ref{set-intersection} to rewrite this probability, showing that
the processes $\overline\YYY_Q$ and~$\overline\YYY_{Q'}$ can be taken
to be nearly independent, and using the fact that, for~$s$ sufficiently
large,
\[
   \pr[\YYY_{P_0}(t/2) \cap \overline\YYY_Q(t/2) \neq \emptyset \giv \ff_s]
    \ \approx\ \pr[\YYY_{P_0}(t/2) \cap \overline\YYY_Q(t/2) \neq \emptyset \giv \ff_{t/2}],
\]
it can be shown that
\[
   \ex\{ [V_{P_0}(t)/L]^2 \giv \ff_s\} \ \sim\ [\ex\{ [V_{P_0}(t)/L] \giv \ff_s\}]^2.
\]
Hence the conditional variance of~$V_{P_0}(t)/L$, given the information in~$\ff_s$, is
small, and thus the value of~$V_{P_0}(t)/L$ is (almost) fixed.  An analogous
argument is used, for instance, in Ball, Sirl \& Trapman~\cite{BallSirlTrap},
where they show that, in an epidemic in a population of large size~$N$,
the proportion of individuals ever infected is
close either to zero or to a non-random value in $(0,1)$.
A by-product of our argument is to identify the solution~$h$ to a
particular integral equation, that appears in Aldous~\cite{Aldous} and also plays a substantial
part in the formula given by Chatterjee \& Durrett~\cite{ChatDurr}, in terms of the Laplace
transform of the branching process limit random variable~$W$; their function~$h$ is just a
time translation of~$h_2$.

The paper is organized as follows. The necessary properties of the
branching processes that approximate the
early stages of the gossip and small worlds processes are established in
Section~\ref{branching}.  The law of large numbers is then proved in Sections
\ref{outline}--\ref{calculations}. The time until~$C$ is completely
covered is investigated in Section~\ref{coverage}, and the paper concludes
in Section~\ref{manifolds}
by extending the results to finite subsets of homogeneous manifolds, such as
rectangles in~$\re^2$.

\section{The branching phase}\label{branching}
 \setcounter{equation}{0}

As in~\cite{BR}, we base our analysis of the coverage process
on the pure growth Markov branching process~$\MBP$, which has neighbourhoods with
centres independently and uniformly positioned in~$C$.
\adb{In this section,} to describe the behaviour of such processes, we
specialize to the case of `flat' manifolds, such as tori, in which
\eq\label{special}
    v_s(\KK) \Eq s^d v(\KK), \qquad s \ge 0.
\en
We later show that this condition can be relaxed substantially, by bounding the
branching processes for more general manifolds
between processes satisfying condition~\Ref{special} that are close enough
for our purposes.

We begin by defining
$M_0(t) := 1 +  \max\{j \ge 0\colon\,\t^*_j \le t\}$ to be the number of islands in
the branching process up to time~$t$, and
\eq\label{M-defs}
   M_l(t) \Eq \sum_{j=1}^{M_0(t)} (t-\t_{j-1}^*)^l
\en
to be the sum of the $l$'th powers of their `radii'.
The evolution of the process is then governed by the differential equations
\eq\label{M-eqns}
   \begin{array}{rl}
   \dfrac{d}{dt} M_1(t) &=\ M_0(t) \quad \mbox{for a.e. }t,  \\[2ex]
   \dfrac{d}{dt} M_i(t) &=\ i M_{i-1}(t), \quad i\ge 2,
   \end{array}
\en
together with a specification of~$M_0$.
Letting~$Z$ denote a unit rate Poisson process,
a small world process is obtained by setting
\eq\label{incidence-area}
   M_0(t) \Eq M_0(0) + Z\left(\r  v(\KK) d\int_0^t M_{d-1}(u)\,du \right)
     \Eq M_0(0) + Z(\r v(\KK) [M_d(t) - M_d(0)]);
\en
for a gossip process, we set
\eq\label{incidence-vol}
   M_0(t) \Eq  M_0(0) + Z\left(\r v(\KK) \int_0^t M_d(u)\,du \right)
                 \Eq  M_0(0) + Z(\r (d+1)^{-1}v(\KK) [M_{d+1}(t) - M_{d+1}(0)]).
\en
In either case, the intensity~$\r$ may depend on~$L$.

Equations \eqref{M-eqns} - \eqref{incidence-vol}
can be rewritten in clearer form by
defining $H_i(t) := M_i(t)\l^i/i!$, for~$\l$ to be suitably chosen, in
which case~\Ref{M-eqns} reduces to
\eq\label{H-eqns}
  \begin{array}{rl}
   \dfrac{d}{dt} H_1(t) &=\ \l H_0(t) \quad \mbox{for a.e. }t, \\[2ex]
   \dfrac{d}{dt} H_i(t) &=\ \l H_{i-1}(t), \quad i\ge 2;
  \end{array}
\en
for the small world process, we have
\eq\label{kappa-area}
    H_0(t) \Eq M_0(0) + Z(d!\r v(\KK) \l^{-d} [H_d(t) - H_d(0)])
      \Eq  M_0(0) +  Z(H_d(t)-H_d(0)),
\en
if $\l = \l_0 := (d!\r v(\KK))^{1/d}$,
and, for the gossip process, we have
\eq\label{kappa-vol}
    H_0(t) \Eq  M_0(0) + Z(d!\r v(\KK) \l^{-d-1} [H_{d+1}(t)-H_{d+1}(0)])
       \Eq  M_0(0) + Z(H_{d+1}(t)-H_{d+1}(0)),
\en
if $\l = \l_0 := (d!\r v(\KK))^{1/(d+1)}$. Note that, since~$\r$ may depend
on~$L$, so may~$\l_0$.

\medskip
\begin{rem}\label{R1}
\adb{The time-scaled process} $\tH(u) := H(u/\l)$ actually satisfies
\eq\label{time-scaled}
   \frac{d}{dt} \tH_i(t) \Eq \tH_{i-1}(t), \quad i\ge 1;\qquad
     \tH_0(t) \Eq  M_0(0) + Z(\tH_{r(d)}(t)-\tH_{r(d)}(0)),
\en
where $r(d) = d$ for the small world and $d+1$ for the gossip process.
Thus, apart from a time change, the processes are the same
for all~$\l$.  Despite this, we retain~$\l_0$ in the subsequent discussion,
in order to emphasize the connection with the original process.
\end{rem}

\medskip
In either case, the equations for~$H = (H_1,H_2,\ldots,H_r)^T$ are of the form
\eq\label{vector-de}
   \frac{dH}{dt} \Eq \l_0 C_r (H + \hat h \e^r) \Eq \l_0 C_r[I + (\hat h/H_r)\e^r(\e^r)^T]\,H,
\en
where $r = r(d)$,  $\e^i$ denotes the $i$-th coordinate vector,
$C_r$ is the $r$-dimensional cyclic permutation matrix satisfying $C_r \e^i = \e^{i-1}$,
$2\le i\le r$, and $C_r \e^1 = \e^r$, and
\eq\label{h-hat-def}
   \hat h(t) \Def  H_0(t) - H_r(t) \Eq Z(H_r(t)-H_r(0)) - H_r(t) + M_0(0) \,.
\en
Without the perturbation~$\hat h$, $H$ would have asymptotically
exponential growth at rate~$\l_0$, and the ratios of its components
would all tend to unity, since the dominant eigenvalue~$1$ of~$C_r$
corresponds to the right eigenvector~$\bone$.  For the arguments
to come, it will be important to show that, with high enough probability,
the asymptotic effect of~$\hat h$ is just to multiply~$H$ by some random
constant --- a branching random variable~$W$ --- which is not too
big.  Unless otherwise specified, we henceforth take
$M_0(0) = 1$ and $M_l(0) = 0$ for all $l\ge1$, so that we start with just
one point~$P_0$ at $t=0$.

\subsection{Growth bounds for the branching process}\label{growth}
Using the maximum norm $\|\cdot\|$ for $r(d)$-vectors, it follows immediately
from~\Ref{vector-de} that
\[
  \frac{d}{dt}\|H(t)\| \Le \l_0 u(t) \|H(t)\|,
\]
with $u(t) := \{1 + (\hat h(t)/H_{r(d)}(t))_+\}$, so that, by a Gronwall argument,
\eq\label{Gronwall}
   \|H(t)\| \Le \|H(t_0)\| \exp\Blb \l_0 \int_{t_0}^t u(v)\,dv \Brb,
\en
for any $0 \le t_0 \le t$.   Thus, in order to bound the
growth of~$H$, we shall need to control the quantity
$\hat h(t)/H_{r(d)}(t)$, which is itself a function of the Poisson process~$Z$.
To do so, we begin with the following lemma, which controls the extreme fluctuations
of~$Z$.

\begin{lem}\label{PP-bound}
   Let~$Z$ be a unit rate Poisson process.  Then we have the following bounds,
uniformly in $t\ge1$:
\eqs
  (1)&& \pr\left[ \sup_{u\ge t} u^{-1}Z(u) \ge 2 \right] \Le c_1\, e^{-t/14}
                  ; \\
  (2)&& \pr\left[ \sup_{u\ge t} u^{1/3}|u^{-1}Z(u)-1| \ge  4 \right] \Le c_2\, e^{-t^{1/3}/5}
                   ; \\
  (3)&& \pr\left[ \inf_{u\ge t} u^{-1}Z(u) \le 1/2 \right] \Le c_3 \, e^{-t/44} .
\ens
Furthermore, for any $U\ge1$ and $0 < \h \le 1/3$ such that $U(2^\h-1) \ge 42\log 2$,
\eqs
   (4)&& \pr\left[ \sup_{u \ge 0} (u\vee1)^{-\shalf(1+\h)}|Z(u)-u| \ge  U \right]
    \Le c_4\, e^{-U/28},
\ens
for a constant~$c_4$.
\end{lem}

\proof
For any $t,\e>0$, set $u_j := t(1+\e)^j$, $j\ge0$.  Then it is immediate that
\[
   \sup_{u_j \le u\le u_{j+1}} u^{-1}Z(u)  \Le Z(u_{j+1})/u_j.
\]
Hence, and by the Chernoff inequalities (\cite{ChungLu}, Theorem~2.3),
\eqs
   \pr\Bigl[ \sup_{u_j \le u\le u_{j+1}} u^{-1}Z(u) \ge 1+2\e \Bigr]
     &\le& \pr[u_j^{-1}Z(u_j(1+\e)) \ge 1 + 2\e] \\
    &\le& \exp\{-\e^2 u_j/(2+3\e)\} \\
    &=& \exp\{-\e^2 u_{j-1}/(2+3\e)\}\,\exp\{-\e^3u_{j-1}/(2+3\e)\},
\ens
and
\eqs
   \pr\Bigl[ \inf_{u_j \le u\le u_{j+1}} u^{-1}Z(u) \le 1 - 2\e \Bigr]
     &\le& \pr\bigl[\{(1+\e)u_j\}^{-1}Z(u_j) \le 1 - 2\e\bigr] \\
    &\le& \exp\{-\e^2 u_j/2\} \\
    &=& \exp\{-\e^2 u_{j-1}/2\}\,\exp\{-\e^3u_{j-1}/2\}.
\ens
Adding over $j\ge0$, the sum is dominated by a geometric progression with
common ratio $\exp\{-\e^3t/(2+3\e)\}$, and so
it follows that
\eq\label{Poisson-bnd}
  \pr\left[ \sup_{u\ge t} |u^{-1}Z(u) - 1| \ge  2\e \right] \Le
     C(\e,t)\exp\{-\e^2 t/(2+3\e)\},
\en
with $C(\e,t) := 2/\{1 - e^{-\e^3 t/(2+3\e)}\}$.  Taking $\e=1/2$ gives the
first inequality, with $c_1 := C(\half,1)$; taking
$\e = 1/4$ gives the third, with $c_3 = C(\quarter,1)$. For the second,
with $t\ge1$, $\e = t^{-1/3}$ gives, in particular,
\[
    \pr\left[ \sup_{t \le u \le 8t} u^{1/3}|u^{-1}Z(u) - 1| \ge  4 \right] \Le
     C(1,1)\exp\{-t^{1/3}/5\},
\]
and thus
\eqa
   \pr\left[ \sup_{u \ge t} u^{1/3}|u^{-1}Z(u) - 1| \ge  4 \right] &\le&
     C(1,1)\sum_{j\ge0}\exp\{-2^jt^{1/3}/5\} \label{new-cross} \\
   &\le& c_2\exp\{-t^{1/3}/5\}, \non
\ena
with $c_2 := C(1,1)/(1 - e^{-1/5})$, since the ratio of successive terms
in~\Ref{new-cross} is at most $e^{-1/5}$.

The fourth inequality is a little trickier.
Taking $t \ge 1$ and $\e = \e_t = U/\{2(2t)^{\shalf(1-\h)}\}$ in~\Ref{Poisson-bnd}, we have
\eqs
    \pr\left[ \sup_{t \le u \le 2t} u^{\shalf(1-\h)}|u^{-1}Z(u) - 1| \ge U \right] &\le&
     C(\e_t,t)\exp\{-2^{-3+\h}U^2t^{\h}/(2 + 3U/2)\} \\
     &\le& C(\e_t,t) \exp\{-Ut^\h/28\},
\ens
since, for $U\ge1$, $U/(2 + 3U/2) \ge 2/7$.
For this choice of~$\e$, $\e^2t$ increases with~$t$, but $\e^3t$ decreases; however,
since $1-e^{-x} \ge (1-e^{-1})\min\{1,x\}$ in $x\ge0$, we have
\[
    C(\e,t) \Le \frac2{1-e^{-1}}\max\Bigl\{1,\frac{2+3\e}{\e^3t}\Bigr\} \Le
       C'(\e,t) \Def \frac2{1-e^{-1}}\max\Bigl\{1,\frac8{\e^3t}\Bigr\},
\]
uniformly in $t\ge1$; the final inequality is immediate if $\e \le 2$,
and, for $\e > 2$, $2+3\e < \e^3$.  Set $q(t) := C'(\e_t,t)\exp\{-Ut^\h/28\}$.  Then, in the sum
$\sjo q(2^j)$, the ratios of successive terms are at most
\[
   2^{(1-3\h)/2} \exp\{-U(2^\h-1)/28\} \Le \sqrt2 \exp\{-U(2^\h-1)/28\} \Le 1/2,
\]
by assumption, so that $\sjo q(2^j) \le 2q(1) \le 2^{11/2}\,\frac{16}{1-e^{-1}}\,\exp\{-U/28\}$.
Since, by an exponential moment inequality, $\pr[Z(1) > U] \le ce^{-U/14}$
with $c := e^{x-1}$ and $x=e^{1/14}$, the proof of the fourth
inequality is complete.
\ep

\bigskip
Based on this lemma, we can now prove growth bounds for the Markov branching process. Here,
we allow for quite general initial conditions.  For ease of reference,
for any $K\ge1$ and $0 < \h < 1$, we define the events
\eqa
   A\ui_{K,s} &:=& \{e^{-\l_0s}\|H(s)\| \le K\};
     \qquad A\ut_{K,s} \Def \{H_{r(d)}(s) \ge K\};\label{ALs-def} \\
   A\uh_{K,\h,s} &:=& \Bigl\{ \sup_{0\le u \le H_{r(d)}(s)} (u\vee1)^{-\shalf(1+\h)}|Z(u)-u|
           \le K^{\shalf(1-\h)} \Bigr\};\label{ALh-def} \\
  A'(K,s) &:=& \{ \exp\{-\l_0 (t-s)(1 + \e_K)\}\|H(t)\| \Le  \|H(s)\|
            \ \mbox{\rm for all}\ t>s \},  \label{ALC-def}
\ena
where $\e_K := 5K^{-1/3}$, and we write
\eq\label{ALs-def-cap}
   A_{K,s} \Def A\ui_{\th K,s} \cap A\ut_{K,s} \cap A\uh_{K,\e_K,s} \ \in\ \ff_s,
\en
where~$\ff_s$ denotes the history of~$\MBP$ up to time~$s$, and~$\th := C_a e^{1/80}$,
with~$C_a$ as defined below.

\begin{thm}\label{H-at-t_L}
For any $K\ge1$ and any $0 \le s < t$, we have
\eqa
 (1) && \pr[\exp\{-\l_0 (t-s)(1 + \e_K)\}\|H(t)\| \le \max\{C_aK, \|H(s)\|\} \ \mbox{\rm for all}\ t>s]
        \non\\
     &&\qquad \ge\ 1 - c_ae^{-K^{1/3}/5}, \non 
\ena
for suitable constants $C_a := 3\exp\{(r(d)!)^{1/r(d)}\}$ and~$c_a$.
Furthermore,
\eqa  
 (2) && \pr[A'(K,s) \giv \ff_s \cap A\ut_{K,s}] \ \ge\ 1 - c_2e^{-K^{1/3}/5}.
     \phantom{HHHHHHHHHHHHHh} \non
\ena
\end{thm}

\proof
  For $K\ge1$, define $\t_{K,s} := \inf\{t\ge s\colon\, H_{r(d)}(t)  \ge K\}$,
and suppose first that $\t_{K,s} > s$.
Then $H_{r(d)}(t) \le K$ if $s \le t \le \t_{K,s}$, and thus, for all such~$t$,
$$
   H_0(t) \Le 1 + Z(K) \Le 3K,
$$
by Lemma~\ref{PP-bound}(1), on a set~$A_1(K)$ of probability at least
$1 - c_1 \exp\{-K/14\}$.
Moreover, from the definition of $H_i(t)$ and by H\"older's inequality
applied to~$M_i(t)$, we have
\eq\label{Hoelder}
   H_i(t) \Le \frac{(r(d)!)^{i/r(d)}}{i!}H_{r(d)}(t)^{i/r(d)} H_0(t)^{1-i/r(d)},
   \qquad 1\le i < r(d).
\en
Hence it follows that, on~$A_1(K)$,
\eq\label{H0}
    \|H(t)\| \Le 3\exp\{(r(d)!)^{1/r(d)}\}K \qquad \mbox{for all}\ s\le t\le \t_{K,s}.
\en

Now, from Lemma~\ref{PP-bound}(2), it follows that
\eqs
  u(t) \Eq 1 + (\hat h(t)/H_{r(d)}(t))_+ &\le& \max\{1,K^{-1} + Z(H_{r(d)}(t))/H_{r(d)}(t)\} \\
    &\le& 1 + 5K^{-1/3} \Eq 1 + \e_K
\ens
for all $t \ge \t_{K,s}$,  on an event~$A_2(K)$ of probability at least $1 - c_2 e^{-K^{1/3}/5}$.
By~\Ref{Gronwall}, this implies that, on~$A_2(K)\cap A_1(K)$, for all $t \ge \t_{K,s}$,
\eq\label{H-bnd-a}
  \|H(t)\| \Le \|H(\t_{K,s})\|\, \exp\{\l_0 (t-s) (1 + \e_K)\}.
\en
If $\t_{K,s} > s$, by \Ref{H0}, this in turn implies that, on~$A_2(K)\cap A_1(K)$,
\eq\label{H-bnd}
  \|H(t)\| \Le 3\exp\{(r(d)!)^{1/r(d)}\}K \exp\{\l_0 (t-s) (1 + \e_K)\}
\en
for all $t \ge s$;
if $\t_{K,s}=s$, we simply have $\|H(t)\| \le \|H(s)\|\, \exp\{\l_0 (t-s) (1 + \e_K)\}$.
This establishes Part~1.

For Part~2, if $H_{r(d)}(s) \ge K$, it follows as above that
\[
   \pr\bigl[u(t) \le 1 + \e_K \ \mbox{for all}\ t\ge s \giv \ff_s\bigr]
      \ \ge\ 1 - c_2e^{-K^{1/3}/5},
\]
and, if this is the case, then
\[
   \exp\{-\l_0 (t-s)(1 + \e_K)\}\|H(t)\| \Le  \|H(s)\| \quad \mbox{for all}\ t>s
\]
follows from~\Ref{Gronwall}.
\ep

 Theorem \ref{H-at-t_L} translates into bounds on the values
of $\hat h(t)$, defined in~\Ref{h-hat-def}.

\begin{cor}\label{epsilon-bnds}
  Given any $\e > 0$, there exists a random variable~$H^\e$ such that
\eqs 
  (1)&& |\hat h(t)| \Le H^\e \exp\{\half\l_0(1+\e)t\}\quad \mbox{a.s.\ for all}\ t > 0.
   \phantom{HHHHHHHHHHi}
\ens
In addition, for any $K\ge1$,
\eqs 
  (2)&& \ex\{|\hat h(t)| I[A'(K,s)] \giv \ff_s\cap A_{K,s}\} \Le
      2\{(\th K)^{1/2} + K\} \exp\{\half\l_0(1+\e_K)t\}.
\ens
\end{cor}

\proof
Note that, from Theorem~\ref{H-at-t_L}\,(1), given any $K_0 \ge 1$, there exists
a.s.\ a (random)~$K \ge K_0$ such that
\[
     \sup_{t>0} e^{-\l_0t(1+\e_K)}\|H(t)\| \ <\ \infty.
\]
Hence, for any $\e>0$,
\eq\label{H-lim-sup}
   H'_\e \Def \max\Bigl\{1,\sup_{t>0} \{e^{-\l_0(1+\e)t}\|H(t)\|\}\Bigr\} \ <\ \infty \  \mbox{a.s.}\,.
\en
This in turn implies that, given any $\e > 0$,
\[
  |\hat h(t) - 1| \Eq |Z(H_{r(d)}(t)) - H_{r(d)}(t)|
     \Le \sup_{0\le u\le H'_{\e/2} \exp\{\l_0(1+\e/2)t\}} |Z(u)-u|\quad\mbox{for all}\ t > 0.
\]
Now, from the law of the iterated logarithm for the
Poisson process, for any $\h > 0$,
\[
    H''_\h \Def \sup_{u\ge1} \{u^{-\h-1/2}|Z(u)-u|\} \ <\ \infty\quad\mbox{a.s.}
\]
Part~1 thus follows because
\[
    \sup_{t\ge1} \Bigl\{e^{-\half\l_0(1+\e)t}|\hat h(t)| \Bigr\}
     \Le 1 + H''_\h H'_{\e/2}\sup_{t\ge1}\{ \exp\{[(\half+\h)(1+\e/2) - \half(1+\e)]\l_0t\} \}
    \ <\ \infty,
\]
if $\h$ is chosen such that $\h(1+\e/2) < \e/4$.

For Part~2, if~$\tZ$ is a Poisson process of rate 1, the Doob--Kolmogorov inequality gives
\[
    \pr\Bigl[ \sup_{0\le u\le U}|\tZ(u)-u| > x\Bigr] \Le \min\{Ux^{-2}, 1\},
\]
and hence $\ex\bigl\{ \sup_{0\le u\le U}|\tZ(u)-u| \bigr\} \le 2\sqrt U$.
Thus, from~\Ref{ALC-def}, we have
\eqa
  \lefteqn{\ex\{|\hat h(t) - \hat h(s)| I[A'(K,s)] \giv \ff_s\cap A_{K,s}\}}\non\\
   &\le&  \ex\Blb \sup_{0\le u\le \|H(s)\| \exp\{\l_0(1+\e_K)(t-s)\}} |\tZ(u)-u|
          \Giv \ff_s\cap A_{K,s} \Brb \non\\
   &\le&  2(\th K)^{1/2} \exp\{\half\l_0(1+\e_K)t\},\label{h-hat-diff}
\ena
where
$\tZ(u) := Z(u + H_{r(d)}(s)) - Z(H_{r(d)}(s))$.  Since also, on $A_{K,s}$,
\[
   |\hat h(s) - 1| \Le H_{r(d)}(s)^{\shalf(1+\e_K)} K^{\shalf(1-\e_K)}
      \Le K\exp\{\half\l_0s(1+\e_K)\},
\]
the proof is completed.
\ep

Recalling~\Ref{vector-de}, and writing $W_*(t) := e^{-\l_0 t}\bone^T H(t)$, we have
\eq\label{W(t)}
    \frac{dW_*}{dt} \Eq \l_0 e^{-\l_0 t} \hat h(t),
\en
and, in view of Corollary~\ref{epsilon-bnds}\,(1), it follows that
$W_*(\infty) := \lim_{t\to\infty} W_*(t)$ exists and is finite a.s.
In fact, the branching process exhibits a non-negative martingale~$\tW$, given by
\eqs
   \tW(t) &:=& e^{-\l_0 t}\sum_{j=0}^{r(d)-1} H_j(t) \Eq W_*(t)
         + e^{-\l_0 t}\{H_0(t) - H_{r(d)}(t)\}
           \Eq  W_*(t) + e^{-\l_0 t}\hat h(t),
\ens
this last by~\Ref{h-hat-def},
for which~$\lti\tW(t) = W_*(\infty)$ also, because of Corollary~\ref{epsilon-bnds}\,(1).
Note that,  in similar fashion, \Ref{W(t)} can also be written as
\eq\label{alternative-de}
    \frac{dW_*}{dt}  \Eq \l_0 e^{-\l_0 t} \{H_0(t) - H_{r(d)}(t)\},
\en
from which, by partial integration, it follows that
$$
   W_*(\infty) \Eq 1 + \int_0^\infty e^{-\l_0 t} \{dH_0(t) - \l_0 H_{r(d)-1}(t)\,dt\},
$$
identifying $W_*(\infty)$ as $r(d) W$, where~$W$ is the limiting random variable
defined in \cite{BR}, Theorem~4.1.

Similarly, integrating~\Ref{vector-de}, and noting that $H(0) = 0$ because
$M_l(0) = 0$ for all $l\ge1$, we have
\[
   e^{-\l_0 C_{r(d)}t}H(t) 
              \Eq \int_0^t \l_0 C_{r(d)} e^{-\l_0 C_{r(d)}u} \e^{r(d)} \hat h(u)\,du,
\]
from which it follows that
\eq\label{H(t)}
   e^{-\l_0 t}H(t) \Eq 
      e^{-\l_0 t}\int_0^t \l_0 C_{r(d)} e^{\l_0 C_{r(d)}(t-u)} \e^{r(d)} \hat h(u)\,du.
\en
Since $\ex \hat h(u) = 1$ for all $u\ge0$, it follows directly that
\eq\label{expectation-bnd}
    e^{-\l_0 t}\| \ex H(t)\| \Le
      e^{-\l_0 t}\int_0^t \l_0 e^{\l_0(t-u)} \|\e^{r(d)}\| \,du \Le 1,
\en
so that $e^{-\l_0 t}H(t)$ is uniformly bounded in expectation.  With some further
effort,  the long term behaviour of the vector $e^{-\l_0 t}H(t)$ can be related in
detailed fashion to that of~$W_*(t)$.

\begin{thm}\label{vector-control}
  As $t\to\infty$, $W(t) := r(d) e^{-\l_0 t}H(t) \to W_*(\infty) \bone$ {\rm a.s.}
Furthermore, for any $0 < s < t$,
\eqs  
\phantom{HH}&& (1)\quad\|W(t) - W(s)\| \Eq O(\exp\{-\half(1-\e)\l_0 s\}) \ \mbox{\rm a.s. }
    \phantom{HHHHHHHHHHHHHH}
\ens
for any $\e > 0$, if $r(d) \le 6$, and
\eqs
\phantom{HH} && (2)\quad\|W(t) - W(s)\| \Eq O(\exp\{-(1-\cos(2\p/r(d)))\l_0 s\}) \ \mbox{\rm a.s. }
   \phantom{HHHHHHHHH}
\ens
if $r(d) \ge 7$.  Finally, for $A'(K,s)$ as defined in~\Ref{ALC-def},
\eqs 
\phantom{HH} && (3)\quad \ex\{\|W(t) - W(s)\| I[A'(K,s)]\giv\ff_s\cap A_{K,s}\}
     \Eq O(K\exp\{-\b_{r(d)}\l_0 s\}),
     \phantom{HHHHHH}
\ens
where $\b_r = \half(1-\e_K)$ if $r \le 6$, and $\b_r = (1-\cos(2\p/r))$ if
$r \ge 7$ and $K^{1/3} \ge 5/\{2\cos(2\p/7)-1\}$.
\end{thm}

\proof
Since the eigenvalues of~$C_r$ are the $r$-th roots $\o_j := e^{2\p (j-1)i/r}$ of unity,
$1\le j\le r$, and the corresponding eigenvectors $e\uj :=
r^{-1/2}(\o_j^{-1},\o_j^{-2},\ldots,\o_j^{-(r-1)},1)$ are (complex) orthonormal,
we can write $e^{-\l_0 t}H(t) = \sum_{j=1}^{r(d)} f_j(t) e\uj$ with
$f_j(t) := e^{-\l_0 t}\{\bar e\uj\}^T H(t)$, and use~\Ref{H(t)} to compute the~$f_j(t)$.

Taking $j=1$ first, pre-multiplying
$e^{-\l_0t}H(t)$  by $r(d)^{1/2}\{e\ui\}^T = \bone^T$ gives
\eqs
    r(d)^{1/2} f_1(t) &=& 
       e^{-\l_0 t}\int_0^t \l_0  e^{\l_0 (t-u)} \bone^T\e^{r(d)} \hat h(u)\,du 
        \Eq \int_0^t \l_0  e^{-\l_0 u} \hat h(u)\,du\, .
\ens
This last expression converges a.s.\ to $W_*(\infty)$ in view of~\Ref{W(t)},
and hence $f_1(t) \to r(d)^{-1/2} W_*(\infty)$ a.s.\ as $t\to\infty$.    Indeed, we have a
little more: for any $s < t$, from Corollary~\ref{epsilon-bnds}\,(1),
\eq\label{H-difference-1}
    r(d)^{1/2} |f_1(t) - f_1(s)|
      \Le \int_s^t \l_0  e^{-\l_0 u} |\hat h(u)|\,du
         \Le \frac{2H^\e}{1-\e} \exp\{-\half(1-\e)\l_0s\}\quad \mbox{a.s. }\,,
\en
for any $\e > 0$.  For $2\le j\le r(d)$, we have
\eqs
  r(d)^{1/2} f_j(t) &=& 
      e^{-\l_0t}\int_0^t \l_0\o_j e^{\l_0\o_j(t-u)} \hat h(u)\,du\,,
\ens
giving
\eqs
  r(d)^{1/2} |f_j(t)| \Le 
      e^{-\l_0t}\int_0^t \l_0 e^{\l_0\r_j(t-u)} |\hat h(u)|\,du\,,
\ens
where $\r_j := \cos\{2\p (j-1)/r_{(d)}\}$.

Corollary~\ref{epsilon-bnds}\,(1) now implies that a.s.
\[
    \int_0^t e^{-\l_0\r u} |\hat h(u)|\,du \Eq
      \begin{cases}
           O(1)  &\mbox{if}\ \r > 1/2;\\
           O\Bigl(\exp\{\l_0t(\half - \r + \e)\}\Bigr)    &\mbox{if}\ \r \le 1/2 \ \mbox{and}\ \e>0.
      \end{cases}
\]
The estimates given in Parts 1 and~2 follow directly, since $\r_2 > 1/2$ if $r(d) \ge 7$,
and
\eq\label{H-limit}
    \lti e^{-\l_0 t}H(t) \Eq r(d)^{-1}W_*(\infty) \bone\quad \mbox{a.s. }
\en
is immediate. The bound given in Part~3 follows in a similar
way, but using Corollary~\ref{epsilon-bnds}\,(2) in place of Corollary~\ref{epsilon-bnds}\,(1);
the final condition on~$K$ is merely to ensure that $\b_7 < \half(1-\e_K)$.
\ep

\begin{rem}\label{W-rk}
Note, in particular, that the distribution of~$W_*(\infty)$ is the same,
for a given value of~$r(d)$,
irrespective of the value of~$\l_0$, since $W(u/\l_0) = r(d) e^{-u}H(u/\l_0)
= r(d) e^{-u}\tH(u)$, with $\tH$ as defined in Remark~\ref{R1}.  We shall
denote a random variable with this distribution by~$W_*^{[r(d)]}(\infty)$,
if the dimension needs to be emphasized.
\end{rem}

For use in Section~\ref{deterministic}, we define
\eq\label{phi-def}
    \f(\th) \Def \ex\{e^{-\th W_*(\infty)}\giv M_0(0)=1,M_j(0)=0,\,j\ge1\}.
\en
The function~$\f$, being the Laplace transform of a branching process limit random variable, can
as usual be expressed as the solution of an implicit equation.  This is based on the observation
that, because of the branching property,
\[
   W_*(\infty) \ =_d\ \sji e^{-\l_0\t_j} W_*^{(j)},
\]
where $(W_*^{(j)},\,j\ge1)$ are independent copies of~$W_*(\infty)$, and $(\t_j,\,j\ge1)$
are the event times in a Poisson process on~$\re_+$.  For the gossip process in dimension~$d$,
the Poisson process has intensity
\[
    \r v(\KK) u^d \Eq \l_0^{d+1}u^d/d!\,,\qquad u>0.
\]
This implies that the Laplace transform~$\f$ satisfies
\eqa
  \f(\th) &=& \exp\Blb - \int_0^\infty \frac1{d!}\{1 - \f(\th e^{-\l_0 u})\}\l_0^{d+1}u^d\,du \Brb
               \non \\
     &=&  \exp\Blb - \int_0^\infty \frac1{d!}\{1 - \f(\th e^{-x})\}x^d\,dx \Brb,
      \label{phi-eqn}
\ena
with $-\f'(0) = \ex\{W_*(\infty)\}\giv M_0(0)=1,M_j(0)=0,\,j\ge1\} = 1$.
If we define $h(t) := 1 - \f(e^t)$, then~$h$ satisfies the equation
\eq\label{h-eqn}
   h(t) \Eq 1 - \exp\Blb - \int_0^\infty \frac{x^d}{d!}h(t-x)\,dx \Brb,
\en
with $\lim_{t\to-\infty} e^{-t}h'(t) = 1$; in the case $d=2$, \Ref{h-eqn}
is the equation that appears in Chatterjee \& Durrett~(\cite{ChatDurr}, Lemma~1.1).  Note that
the functions $\f= \f_d$ and~$h = h_d$ thus depend only on the dimension~$d$, and not
on the choice of neighbourhoods, and hence that this is true also of the distribution
of~$W_*(\infty)$.

It follows from~\Ref{phi-eqn} that the lower tail of~$W_*(\infty)$ can be easily bounded:
\eqa
   \pr[W_*(\infty) \le w] &\le& e\f(1/w)
    \Eq e\exp\Blb - \int_0^\infty \frac1{d!}\{1 - \f(e^{-x}/w)\}x^{d}\,dx \Brb  \non\\
      &\le& e\exp\Blb - \int_0^{-\log w} \frac1{d!}\{1 - \f(1)\}x^{d}\,dx \Brb  \non\\
      &\le& e\Bl\exp\Blb - c(\log(1/w))^{d+1} \Brb \Br, \label{W-lower-tail-bnd}
\ena
for $c := \{1-\f(1)\}/(d+1)!$.  Thus $\pr[W_*(\infty) \le w]$ goes to zero faster than
any power of~$w$ as $w\to0$ for all $d\ge1$.  For the upper tail,
\eq\label{W-upper-tail-bnd}
    \pr[W_*(\infty) \ge w] \Le 1/w
\en
is immediate from Markov's inequality.

Our final result of the section is a lower bound on the growth of~$H$,
needed later to show that $\pr[A\ut_{K,s}]$ 
is large enough, if~$s$ is not too small.

\begin{lem}\label{lower-bnd}
Let~$\t^{r(d)}_K$ denote the time taken for~$H_{r(d)}$ to first reach a given level~$K$.
Then there is a constant~$c_c := 2(r(d)!)^{1/r(d)}/\log(6/5)$ such that
\[
   \pr[\t^{r(d)}_K \ge c_c \l_0^{-1}R] \Le 2Ke^{-R},
\]
for all $R \ge1$.
\end{lem}

\proof
Define $t_0 := \l_0^{-1}(r(d)!)^{1/r(d)}$, so that, starting~$\MBP$ with
a single particle at~$0$, we have $H_{r(d)}(t_0) \ge 1$.  Then it is immediate that
the process~$M_0$ is stochastically bounded below by a process~$\tM_0$,
where $(\tM_0(nt_0),\,n\ge0)$ is a Galton--Watson process with $\tM_0(0)=1$
and with offspring
distribution $p_1=(1-p_2)=e^{-1}$, and that, for $nt_0 \le t < (n+1)t_0$,
$\tM_0(t) = \tM_0(nt_0)$.  Furthermore, for all~$n\ge1$,
$H_{r(d)}(nt_0) \ge \tM_0((n-1)t_0)$. Now, since $\Bi(m,\half)\{[m/2,m]\} \ge 1/2$
for all $m\ge1$, a calculation shows that
\[
   \ex\{1/\tM_0((n+1)t_0) \giv \tM_0(nt_0) = m\} \Le m^{-1}(\half + \third),
\]
from which it follows that $((6/5)^n/\tM_0(nt_0),\,n\ge0)$ is a non-negative
supermartingale, with initial value~$1$.  So, defining $\n_K := \min\{n\colon\,
\tM_0(nt_0) \ge K\}$, and noting also that $\tM_0(\n_Kt_0) \le 2K$, it follows
from the optional stopping theorem that
\[
    \ex\{(6/5)^{\n_K}\} \Le 2K.
\]
It thus follows that $\pr[\n_K \ge R / \log(6/5)] \le 2Ke^{-R}$,
for any $R>0$.  As $\t^{r(d)}_K/t_0 \le \n_K+1$ we have
\[
   \pr[\t^{r(d)}_K \ge t_0\{1 + R / \log(6/5)\}] \le 2Ke^{-R}
\]
also.  Since $R + \log(6/5) \le 2R$ in $R\ge1$, the lemma follows.
\ep

\subsection{Intersection asymptotics}
The branching process~$\MBP$ gives a useful approximation to~$Y$ as long as it has
only few pairs of self intersecting islands.  Thus our asymptotics for the branching
process are of most interest in the time before self-intersections in~$\MBP$ become
plentiful.  To see when this is, we need
to derive formulae for the mean number of pairs of self intersecting
islands of~$\MBP$ at a given time~$t$, and the number of islands of one process~$\MBP_1$
that intersect islands of another, independent process~$\MBP_2$.

In order to do so, we strengthen the assumption that $v_s(\KK) \sim s^d v(\KK)$ as $s\to0$
by assuming that
\eq\label{volumes}
   |s^{-d}v_s(\KK) / v(\KK) - 1| \Le c_g (s\{v(\KK)/L\}^{1/d})^{\g_g},\qquad s > 0,
\en
for some $\g_g > 0$.
Two subsets $\KK(P,t)$ and~$\KK(Q,u)$ intersect when $P\in\KK(Q,t+u)$,
or, equivalently, when $Q\in\KK(P,t+u)$, so that the probability of intersection if $P$
and~$Q$ are chosen independently and uniformly distributed on~$C$ (with respect to the intrinsic volume),
is given by
\eq\label{intersection-prob}
   \ppp_L(t,u) \Def L^{-1}v_{t+u}(\KK) \Eq L^{-1}v(\KK)(t+u)^d(1 + R_L(t,u)),
\en
where
\eq\label{intersection-error}
   |R_L(t,u)| \Le c_g \{(t + u)(v(\KK)/L)^{1/d}\}^{\g_g}.
\en

We  begin by supposing that $c_g=0$.  In this case, the number~$N(t)$ of self intersecting
pairs of islands in a branching process at time~$t$, conditional on having
$M_0(t) = n+1$ islands originating at times $\t^*_0=0,\t^*_1,\ldots,\t^*_n
\le t$, has mean given by
\eqa\label{self-int-mean}
   \ex \{N(t) \giv M_0(t) = n+1, \t^*_1,\ldots,\t^*_n\} &=& \m(n+1;t,t-\t^*_1,\ldots,t-\t^*_n)
      \non\\
   &:=&  \sum_{i=0}^{n-1} \sum_{j=i+1}^n L^{-1}v(\KK)
         \sum_{l=0}^d {d\choose l} (t-\t^*_i)^l (t-\t^*_j)^{d-l}\phantom{XX}\non\\
     &=& \half L^{-1}v(\KK)\sum_{l=0}^d {d\choose l} \{M_l(t) M_{d-l}(t) - M_d(t)\}.
\ena
Similarly, again if $c_g = 0$, for two independent branching processes, one that
has developed to time~$t$ and has $M_0(t)=m+1$ islands originating
at times $0,\t^*_1,\ldots,\t^*_m \le t$, and the other that has developed to time~$s$ and has
$\tM_0(s)=n+1$ islands originating at times $\s^*_0=0,\s^*_1,\ldots,\s^*_n \le s$, the conditional
mean of the number~$N(t,s)$ of intersecting pairs is given by
\eqa
   \lefteqn{\ex \{N(t,s) \giv M_0(t) = m+1, \t^*_1,\ldots,\t^*_m,\,\tM_0(s)=n+1,\s^*_1,\ldots,\s^*_n \}
                 }\non\\  [2ex]
   &&\Eq \m'(m+1,n+1;t,t-\t^*_1,\ldots,t-\t^*_m;s,s-\s^*_1,\ldots,s-\s^*_n) \non \\
   &&\Def \sum_{i=0}^{m} \sum_{j=0}^n L^{-1}v(\KK)
        \sum_{l=0}^d {d\choose l} (t-\t^*_i)^l (s-\s^*_j)^{d-l}   \non\\
    &&\Eq L^{-1}v(\KK)\sum_{l=0}^d {d\choose l} M_l(t) \tM_{d-l}(s), \label{int-mean}
\ena
where $\tM_l(s) := \sum_{j=0}^{n} (s-\s^*_j)^l$.  Thus the quantities~$M_l(\cdot)$
of the previous section are exactly the quantities needed for making such
computations.
If $c_g > 0$, we instead have
\eqa
  && |\ex \{N(t) \giv M_0(t) = n+1, \t^*_1,\ldots,\t^*_n\} - \m(n+1;t,t-\t^*_1,\ldots,t-\t^*_n)|
                  \phantom{XXXXXXX}\non\\
     &&\hskip1in \Le {n\choose2}  c_g \{(2t)^d v(\KK)/L\}^{1 + \g_g/d};  \non\\
  && \Bigl|\ex \{N(t,s) \giv M_0(t) = m+1, \t^*_1,\ldots,\t^*_m,\,\tM_0(s)=n+1,\s^*_1,\ldots,\s^*_n \}
            \phantom{XXXXXXX}     \label{int-means-curved}\\ 
  && \quad\qquad\mbox{} - \m'(m+1,n+1;t,t-\t^*_1,\ldots,t-\t^*_m;s,s-\s^*_1,\ldots,s-\s^*_n)\Bigr| \non\\
    &&\hskip1in \Le mn c_g \{(2t)^d v(\KK)/L\}^{1 + \g_g/d}.  \non
\ena

Now, from Theorem~\ref{vector-control}, we have
\[
    M_l(t) M_{d-l}(t)\ \sim\ l!(d-l)!\l_0^{-d}r(d)^{-2}e^{2\l_0t},
\]
so that the mean number of self-intersections is small up until times~$t$
at which
\eq\label{Lambda-def}
     e^{2\l_0t}\ \asymp\ \L \Def L\l_0^d/v(\KK),
\en
or, equivalently, for
\eq\label{t-lambda-def}
   t \Eq t_{\L,x} \Def \frac1{2\l_0}\{\log\L + x\}.
\en
 As can be seen from~\Ref{vector-de}, $\l_0^{-1}$ is the time scale
for the growth of~$H(\cdot)$, or equivalently, in view of \Ref{kappa-area}
and~\Ref{kappa-vol}, the time scale on which the first new contact occurs.
Correspondingly, $v(\KK)\l_0^{-d}$ represents the scale for the size of the
initial island at the time when the first new contact occurs.
For the sort of asymptotics to be discussed here, it is natural to require
that this size is small when compared to the total size~$L$ of~$C$, so that
many islands are involved in the covering of~$C$;
hence we shall think of the ratio~$\L = L / \{v(\KK)\l_0^{-d}\}$ as being large.

\section{The deterministic phase}\label{deterministic}
 \setcounter{equation}{0}
\subsection{Outline}\label{outline}
In this section, we show that the development of the proportion~$L^{-1}V(t)$ of the volume of~$C$
that is covered at time~$t$ grows more or less deterministically, once the initial
stages have passed,  under the assumption~\Ref{volumes}.
As noted above, up to times of the form
$\half\l_0^{-1}\{\log\L + x\}$, for~$x$ fixed, there are few self-intersections in
the branching process~$\MBP$, so that, for such times, calculations made using the
branching process can be expected to give close to the right answers for the
small world and gossip processes as well.
Differences between $\MBP_{P_0}$ and~$\YYY_{P_0}$ arise because some islands
in~$\MBP_{P_0}$ do not appear in~$\YYY_{P_0}$, and are labelled as `ghosts';
the $j$-th new branching process contact contributes to~$\YYY$ (designated by
$G_j=0$) only if it originates
from a real (non-ghost) island, and starts at a point outside~$\YYY$, and multiple counting
from overlapping islands is also prohibited.

The initial similarity between a small world model and its branching process approximation was
exploited in \cite{BR}, when approximating the distribution
of the distance $d_{SW}(P,P')$ between two randomly chosen points $P$ and~$P'$ in the small
world graph.  The key observation is that $d_{SW}(P,P') > 2t$ exactly when the sets $\YYY_P(t)$
and~$\YYY_{P'}(t)$ are disjoint.  The discussion above indicates that the sets $\YYY_P(t)$
and~$\YYY_{P'}(t)$ can be replaced with little error, for the calculation of probabilities,
by the sets $\YYYs_P(t)$ and~$\YYYs_{P'}(t)$ generated by independent branching processes
$\MBP_P$ and~$\MBP_{P'}$,
as long as~$t$ is of the form $\half\l_0^{-1}\{\log\L + x\}$, which
is enough for our purposes.
However, in small world models, the argument is more complicated,
because an intersection between an island $J$ of $\MBP_P(t)$ and an island~$J'$ of $\MBP_{P'}(t)$
does not correspond to an intersection of $\YYY_P(t)$ with~$\YYY_{P'}(t)$
when either $J \subset J'$ or $J' \subset J$, since such a constellation cannot occur in
the small world process.  Thus, for~$P$ an independent uniform point of~$C$, we have
\eq\label{main-asymp-a}
    \pr_{P_0}[d_{SW}(P_0,P) \le 2t] \Eq \pr[\YYY_{P_0}(t) \cap \YYY_{P}(t) \ne \emptyset]\ \sim\
       \pr[\YYYs_{P_0}(t) \cap^* \YYYs_{P}(t) \ne \emptyset],
\en
for independent $\MBP_{P_0}$ and~$\MBP_P$,
where the notation~$\mbox{}\cap^*\mbox{}$ is used to denote this special mode of intersection.
The complement of the latter probability was approximated,
for small world processes, in \cite{BR}, Theorem~4.2.

The connection with the current problem is that, in the small world model,
\eq\label{main-asymp-b}
  \pr_{P_0}[d_{SW}(P_0,P) \le 2t] \Eq  L^{-1} \ex V_{P_0}(2t),
\en
if~$P$ is independently and uniformly chosen on~$C$,
and, similarly,
\eq\label{main-asymp-c}
  L^{-1} \ex \{V_{P_0}(2t) \giv\ff_s\} \Eq \pr_{P_0}[d_{SW}(P_0,P) \le 2t \giv \ff_s],
\en
where $\ff_s$ denotes the history of~$\YYY_{P_0}$ up to time~$s$.  Then
\eq\label{main-asymp-c1}
   \pr_{P_0}[d_{SW}(P_0,P) \le 2t \giv \ff_s]
     \ \sim\ \pr[\YYYs_{P_0}(t) \cap^* \YYYs_{P}(t) \ne \emptyset \giv \ff_s] \ =:\ \p(s,t).
\en
where the asymptotic equivalence can be expected much as for~\Ref{main-asymp-a}.
The aim is now to demonstrate that the quantity
$L^{-1}V_{P_0}(2t)$ stays close to its expectation $L^{-1} \ex \{V_{P_0}(2t) \giv \ff_s\}$,
showing that the volume process develops in almost
deterministic fashion from~$s$ onwards as long as~$s$ is sufficiently large.
For this, it is enough to show that the conditional
variance $\var\{L^{-1}V_{P_0}(2t) \giv \ff_s\}$ becomes small with~$s$.  Now the expectation
$L^{-1} \ex \{V_{P_0}(2t) \giv \ff_s\}$ is approximated, as above, by~$\p(s,t)$.
For the mean square, we simply note that, for $t>s$,
\eq\label{pair-prob-1}
  L^{-2} \ex \{[V_{P_0}(2t)]^2 \giv \ff_s\}
    \Eq \pr_{P_0}[\{d_{SW}(P_0,P) \le 2t\} \cap^* \{d_{SW}(P_0,P') \le 2t\} \giv \ff_s],
\en
where $P$ and~$P'$ are two independent uniform points of~$C$,
and then, again much as for~\Ref{main-asymp-a}, establish the approximation
\eqa
  \lefteqn{\pr_{P_0}[\{d_{SW}(P_0,P) \le 2t\} \cap \{d_{SW}(P_0,P') \le 2t\} \giv \ff_s] }\\
    &&\ \sim\ \pr[\{\YYYs_{P_0}(t) \cap^* \YYYs_{P}(t) \ne \emptyset\} \cap
      \{\YYYs_{P_0}(t) \cap^* \YYYs_{P'}(t) \ne \emptyset\} \giv \ff^{*}_s],
  \label{pair-prob}
\ena
using three independent processes $\MBP_{P_0}$, $\MBP_{P}$ and~$\MBP_{P'}$,
where $\ff^{*}_s$ denotes the history of~$\MBP_{P_0}$ up to~$s$.  Since, for~$s$
large enough, 
the development of the branching process~$\MBP_{P_0}$
after~$s$ is almost deterministic, the statistics of the set~$\YYYs_{P_0}(t)$ are already
almost determined at time~$s$. Hence the probability in~\Ref{pair-prob} is close to
\[
   \{\pr[\YYYs_{P_0}(t) \cap^* \YYYs_{P}(t) \ne \emptyset \giv \ff^{P_0}_s]\}^2\ =:\ \p^2(s,t),
\]
asymptotically equivalent to the square of the conditional mean.
In consequence, the conditional variance is small, as required.  The remainder of the section
consists of making an analogous argument precise, in the context of gossip processes.
For small world processes, the corresponding estimates can be deduced using the methods
of \cite{BR}, Section~4.

\subsection{Constructions}\label{const}
For the gossip process, the
following construction of~$\YYY$ is useful.  First, for any $(P,t) \in
C \times [0,\infty)$ let $S_{\KK}(P,t) \subset C \times [0,\infty)$ denote the set
whose $C$--section $\s_{C,u}\{S_{\KK}(P,t)\}$ at~$u$ is $\KK(P,u-t)$ if $u\ge t$
and~$\emptyset$ otherwise.
Let~$\Pi$ denote a marked Poisson
process on $C \times [0,\infty)$ with constant intensity~$\r$, and with marks
\iud\ in~$C$.  Take $P_0 \in C$ to be the initial point; set
$S_0 := S_{\KK}(P_0,0)$.  Then define $S_1$ to be the set
\[
   S_1 \Def  S_0   \bigcup_{j\ge1} S_{\KK}(P_{1j},\t_{1j}) ,
\]
where the points of~$\Pi$ in $S_0$ occur at locations and times $(Q_{1j},\t_{1j})$
and have marks~$P_{1j}$, $j\ge1$.  Then recursively, for $l\ge2$, define
\[
   S_{l} \Def  S_{l-1} \bigcup_{j\ge1} S_{\KK}(P_{lj},\t_{lj}) ,
\]
where the points of~$\Pi$ in $S_{l-1} \setminus S_{l-2}$ occur at locations and times
$(Q_{lj},\t_{lj})$ and have marks~$P_{lj}$,  $j\ge1$. Then we can define
\eq\label{Y-def-2}
    \YYY(t) \Def \s_{C,t}\{S_\infty\}, \qquad 0\le t\le T,
\en
where $S_\infty := \cup_{l\ge0}S_l$.  Clearly,
$S_\infty \cap\{C\times[0,T]\} = S_l \cap\{C\times[0,T]\}$ for some $l\ge0$,
since~$\Pi$ has a.s.\ only finitely many points in~$C\times[0,T]$.

Thinking of the points $(Q_{lj},\t_{lj})$ with marks~$P_{lj}$ as the $l$-th
generation descendants of the initial individual~$P_0$, we note that the
path of the process~$\YYY$
can be identified with that of a Markov branching process~$\MBP$ until the first time~$t$
at which there is a non-empty intersection between two of the sets $S_{\KK}(P_{lj},t)$,
for $l\ge0$ and $j\ge1$ such that $\t_{lj} \le t$; this is because of the independence
of Poisson realizations on disjoint sets.   Thus the construction of~$\YYY$ also yields
a part of a coupled~$\MBP$.

As an alternative,
the construction of $\Pi$ and~$\YYY$ can be replicated starting from the branching
process~$\MBP$. Start with~$P_0$, and, writing $S_0 := S_{\KK}(P_0,0)$, assign Poisson
points $(\tQ_{1j},\tit_{1j})$, $j\ge1$, to $S_0 \cap \{C\times[0,T]\}$ with
intensity~$\r$, and with marks \iud\ on~$C$; denote these by
$I_{1j} := (\tQ_{1j},\tit_{1j};\tP_{1j})$, $j\ge1$.  These are the first generation
descendants of the individual at~$P_0$ that are born before~$T$.  Repeat the process
recursively, at the $l$-th step, $l\ge2$, assigning Poisson points and marks to
each~$S_{l-1,j}\cap \{C\times[0,T]\}$, where $S_{l-1,j} := S_{\KK}(\tP_{l-1,j},\tit_{l-1,j})$,
yielding the $l$-th generation descendants
$I_{lj} := (\tQ_{lj},\tit_{lj};\tP_{lj})$, $j\ge1$, of the individual at~$P_0$ that
are born before~$T$; as before, write $S_l := S_{l-1}\cup\{\cup_{j\ge1} S_{lj}\}$.
There are in general more of these descendants than there are Poisson points and marks
in~$S_\infty\cap\{C\times[0,T]\}$,
and the labelling is typically different.  To recover the Poisson process~$\Pi$ and
its marks, first identify the set $\{I_{1j},\, j\ge1\}$ with
$\{( Q_{1j},\t_{1j};P_{1j}),\, j\ge1\}$, the first generation descendants of the individual
at~$P_0$ in the gossip process.  Then label
each~$I_{1j}$ with $G_{1j} \in \{0,1\}$, according as whether the point is to be
treated as real or a ghost; $G_{1j} = 1$ if $(P_{1j},\t_{1j}) \in S_0$, and, if there
is any~$j$ with~$G_{1j}=1$, the union in~\Ref{Y-def-2} is not disjoint beyond~$\t_{1j}$.

For the remaining construction, the descendants~$I_{lj;r}$ of individual~$I_{lj}$ in $[0,T]$
are the  Poisson points and marks from~$S_{lj}\cap \{C\times[0,T]\}$
chosen above. If $G_{lj} = 1$,
label them all with $G_{lj;r}=1$ --- as descendants of ghosts, they are themselves ghosts.
If $G_{lj} = 0$, label those of the~$I_{lj;r}$ for which
\[
   (\tQ_{lj;r},\tit_{lj;r}) \in  S_{l-1} \bigcup_{j'=1}^{j-1} S_{lj'}
\]
with $G_{lj;r} = 1$ also, since they lie in a part of $C \times [0,T]$ that has already
been covered, and therefore do not belong to~$\Pi$.  Do the same if
\[
   (\tP_{lj;r},\tit_{lj;r}) \in  S_{l-1} \bigcup_{j'=1}^{j} S_{lj'},
\]
since such points do not generate new descendants in~$\YYY$, because then
$S_{\KK}(\tP_{lj;r},\tit_{lj;r})$ is contained in the already covered region
of~$C \times [0,\infty)$, even though here the points $(\tQ_{lj;r},\tit_{lj;r})$
do represent points of~$\Pi$.  The remaining points determine the points of~$\Pi$ in
$S_{lj} \setminus \Blb S_{l-1} \bigcup_{j'=1}^{j-1} S_{lj'} \Brb$ that carry the
label $G_{lj;r} = 0$ and are in $C \times [0,T]$.  The set of points
\[
   \Blb (\tQ_{lj;r},\tit_{lj;r}),\,j\ge1, r\ge1:
     (\tQ_{lj;r},\tit_{lj;r}) \notin  S_{l-1} \bigcup_{j'=1}^{j-1} S_{lj'} \Brb,
\]
together with their marks $\tP_{lj;r}$,
can now be identified with the set of points and marks
$\{(Q_{l+1,j},\t_{l+1,j};P_{l+1,j}),\,j\ge1\}$, after a suitable re-indexation
(and with the labels~$G_{l+1,j}$ correspondingly defined), recovering the gossip
process~$Y$.

Note that the same construction can also be used starting from two independent
branching processes $X^{(1)*}$ and~$X^{(2)*}$ with initial points $P\ui_0$
and~$P\ut_0$, provided that, in each generation, the descendants of both
individuals are listed together.  The resulting gossip process describes the
informed regions at each time~$t$, when the information spreads from two
initial sources at $P\ui_0$ and~$P\ut_0$.  Note that the order in which
individuals appear in the list can influence the realization of~$Y$ that
is obtained, so that, for instance, the set $Y_{P\ui_0}(t)$ obtained
from~$X^{(1)*}$ alone may not be identical to the subset of points within
distance~$t$ of~$P\ui_0$ in the set $Y_{P\ui_0,P\ut_0}(t)$ obtained from
$X^{(1)*}$ and~$X^{(2)*}$ together.  However, it is shown in the proofs below
that the differences, which only occur as a result of sets in the branching
processes overlapping, are not significant for the ranges of~$t$ under
consideration here.

For the calculations to come, we next need to show that analogues of the asymptotic
equivalences in \Ref{main-asymp-c1} and~\Ref{pair-prob} hold for
gossip processes.  Here, the argument is a little simpler than for small world
processes.  A point~$P$ has been informed from~$P_0$ at time~$2t$ exactly when the
set~$\YYY_{P_0}(t) \subset C$ of points informed from~$P_0$ by time~$t$ intersects the
set of points~$\bC_{P}(t;2t) \subset C$ from which the information will reach~$P$
by time~$2t$ if it has reached $\bC_{P}(t;2t)$ by time~$t$.  Now the set~$\bC_{P}(t;2t)$
is determined by the points $(\t^t_j,P^t_j)$ of~$\Pi$ in $C\times(t,2t]$, together with their
associated marks~$Q^t_j$, and can be constructed from them in exactly the same way
as~$\YYY_{P_0}(t)$ is constructed from the points $(\t^0_j,P^0_j)$ of~$\Pi$ in $C\times(0,t]$,
together with their associated marks~$Q^0_j$, except that time is run backwards
from $2t$ to~$t$, and the roles of $P^t_j$ and~$Q^t_j$ are swapped.  Thus, and because
the neighbourhoods~$\KK$ were so chosen that $P \in \KK(Q,t)$ exactly when $Q \in \KK(P,t)$,
the set
$\bC_{P}(t;2t)$ has the same distribution as~$\YYY_{P}(t)$, and $\YYY_{P_0}(t)$
and~$\bC_{P}(t;2t)$ are {\it independent\/}.
Thus we can deduce expressions based on which analogues of \Ref{main-asymp-c1} and~\Ref{pair-prob}
can be justified, using the branching process approximations of Section~\ref{branching}:
for any~$s < t$,
\eqa\label{main-asymp-goss}
  L^{-1} \ex \{V_{P_0}(2t)\giv \ff_s\}
      &=& \pr_{P_0}[\YYY_{P_0}(t) \cap \bC_{P}(t;2t)\ne \emptyset\giv \ff_s]; \\
  L^{-2} \ex \{[V_{P_0}(2t)]^2 \giv \ff_s\}
      &=& \pr_{P_0}[\{\YYY_{P_0}(t) \cap \bC_{P}(t;2t)\ne \emptyset\}
         \cap \{\YYY_{P_0}(t) \cap \bC_{P'}(t;2t)\ne \emptyset\} \giv \ff_s], \non
\ena
where $P$ and~$P'$ are independent uniform points of~$C$. The argument is now
primarily aimed at justifying the replacement of $\YYY_{P_0}(t)$, $\bC_{P}(t;2t)$
and~$\bC_{P'}(t;2t)$ in the formulae above by independent copies of~$\YYYs_Q(t)$ with
appropriate choices of~$Q$, so that computations can conveniently be made.

\subsection{Calculations}\label{calculations}
As remarked at the start of Section~\ref{branching}, it is useful to be able to bound~$\MBP$ above
and below by branching processes having $c_g=0$, to which we can then apply the results
of that section.  We do this by constructing processes $\Xm$ and~$\Xp$, for times
$t \ge s$, with the same initial conditions as~$\MBP$ at~$s$, and using the same underlying
Poisson process~$Z$; the time~$s$ and the initial conditions can be freely chosen.
For~$\MBP$, the quantity $v(\KK)M_d(u)$ has to be replaced
by~$V^*(u)$ in the argument of~$Z$ in~\Ref{incidence-vol}, and, by~\Ref{volumes}, we have
\[
    v(\KK)M_d(u)\{1-\h_\L\} \Le V^*(u) \Le v(\KK)M_d(u)\{1+\h_\L\},
\]
where $\L$ is as defined in~\Ref{Lambda-def} and
\eq\label{eta-L-def}
   \h_\L \Def  c_g\Bl \frac{3\log\L}{2\L^{1/d}} \Br^{\g_g},
\en
uniformly in $0 \le u \le (3/2)\l_0^{-1}\log\L$.  Hence we can define~$\Xp$ as in
Section~\ref{branching}
by using $\r^+ := \r\{1+\h_\L\}$ as the contact rate per unit volume, and~$\Xm$
with $\r^- := \r\{1-\h_\L\}$, in which case
\eq\label{Xpm-bnds}
   M_l^-(t) \Le M_l(t) \Le M_l^+(t) \quad\mbox{for all}\ 0\le t\le (3/2)\l_0^{-1}\log\L
    \ \mbox{and}\ 0 \le l \le d+1.
\en
We shall use $\Xm$ and~$\Xp$ extensively to bound quantities associated with~$\MBP$,
and write
\eq\label{lambda-0-pm}
    \l_0^{\pm} \Def \l_0\{1\pm\h_\L\}^{1/d}
\en
for the corresponding growth exponents;
for convenience, we shall assume henceforth that~$\L$ is large enough that $9\h_\L\log\L \le 1$.

We now continue with the following Poisson approximation result (see \cite{BR}, Theorem~3.1),
which can be simply proved using the Stein--Chen method.

\begin{lem}\label{Poisson}
   Let~$n$  $\KK$-islands of radii $t_1,\ldots,t_n \le (3/2)\l_0^{-1}\log\L$ have centres
independently and uniformly distributed on~$C$, and let $N_n$ denote the number of
pairs of them that intersect.  Then
\[
   \dtv(\law(N_n),\Po(\ex N_n)) \Le 4np_\L^+,
\]
where, recalling~\Ref{intersection-prob},
\eq\label{plp-def}
   p_\L^+ \Def \frac{\{3\log\L\}^d}{\L} \Blb 1 
        +  c_g\Bl \frac{3\log\L}{\L^{1/d}} \Br^{\g_g} \Brb
   \Eq \frac{\{3\log\L\}^d}{\L}\{1 + 2^{\g_g}\h_\L\}
\en
is an upper bound for the probability that two independently positioned $\KK$-islands of radius
at most $(3/2)\l_0^{-1}\log\L$ intersect.  Similarly, for two independent
collections of $\KK$-islands, one with $m$ and one with~$n$ islands, having radii $t_1,\ldots,
t_m  \le (3/2)\l_0^{-1}\log\L$ and $u_1,\ldots,u_n  \le (3/2)\l_0^{-1}\log\L$ respectively,
the number~$N_{mn}$ of intersecting pairs satisfies
\[
   \dtv(\law(N_{mn}),\Po(\ex N_{mn})) \Le 2(m+n)p_\L^+.
\]
\end{lem}

Thus, when finding the probability of there being an intersection between the islands of
$\MBP_{P_1}(t)$ and~$\MBP_{P_2}(t)$, where $\MBP_{P_1}$ and~$\MBP_{P_2}$ are independent,
this Poisson approximation offers an approach.  It is exploited in the following result,
in which $M\ui$ and~$M\ut$ denote the quantities~\Ref{M-defs} derived from
$\MBP_{P_1}(t)$ and~$\MBP_{P_2}(t)$.

\begin{lem}\label{br-pr-intersect}
  Define
\[
   P^*[m\ui,m\ut] \Def
      \pr[\YYYs_{P_1}(t) \cap \YYYs_{P_2}(t) \ne \emptyset \giv M\ui(t) = m\ui, M\ut(t) = m\ut].
\]
Then, for $t \le (3/2)\l_0^{-1}\log\L$, we have
\eqs
   &&\Bigl| P^*[m\ui,m\ut] -
        \Bigl( 1 - \exp\Bigl\{-\sjd {d\choose j} m_j\ui m_{d-j}\ut L^{-1}v(\KK)\Bigr\}\Bigr)\Bigr| \\
    &&\qquad  \Le 2(m\ui_0 + m\ut_0)p_\L^+
      +  m\ui_0 m\ut_0 \frac{\{3\log\L\}^d}{\L}\, 2^{\g_g}\h_\L.
\ens
In particular, for $t = t_{\L,x} := \half\l_0^{-1}\{\log\L + x\}$, as in~\Ref{t-lambda-def},
for $x \le \half\log\L$ and for $K\ge1$, this gives
\eqa\label{nice-prob}
   \lefteqn{\Bigl| P^*[M\ui(t),M\ut(t)] -
        \Bigl( 1 - \exp\Bigl\{-\sjd {d\choose j}  M_j\ui(t) M_{d-j}\ut(t) \l_0^d \L^{-1}
            \Bigr\}\Bigr)\Bigr|} \non\\
      &&\qquad\Le C\Blb K\L^{-\shalf + \stq\e_K}\{\log\L\}^{d}e^{x/2} 
            + K^2\L^{-\g_g/d+3\e_K/2}\{\log\L\}^{d+\g_g} \Brb,
\ena
for suitable choice of~$C$, except on an event of probability at most~$2c_ae^{-K^{1/3}/5}$.
\end{lem}

\proof
The first part of the lemma is an immediate consequence of Lemma~\ref{Poisson}, together
with~\Ref{int-means-curved}.  For the second part, we use~\Ref{Xpm-bnds} to bound $M\ui_0(t)$
and~$M\ut_0(t)$ above by $M\uip_0(t)$ and~$M\utp_0(t)$, and then Theorem~\ref{H-at-t_L}\,(1)
to show that, for $l=1,2$,
\[
  M_0^{(l)}(t) \Le C_a K e^{\l_0^+t(1+\e_K)} \Le eC_aK e^{\l_0 t(1+\e_K)}
\]
for all~$0\le t \le  (3/2)\l_0^{-1}\log\L$, except on a set of probability at most
$c_ae^{-K^{1/3}/5}$, since, for~$t$ in this range,
\eq\label{lpt-bnd}
   \l_0^+t(1+\e_K) \Le \{\l_0t + (3/2)\h_\L\log\L\}(1+\e_K) \Le \l_0t(1+\e_K)+1,
\en
if $9\h_\L\log\L \le 1$, since also $\e_K \le 5$ for $K\ge1$.
\ep

We now aim to show that, for $t = t_{\L,x}$ and $s < t$, the
conditional probabilities of actual interest for~\Ref{main-asymp-goss},
\eq\label{p1-def}
   p_1(s,t) \Def \pr_{P_0}[\YYY_{P_0}(t) \cap \bC_{P}(t;2t)\ne \emptyset\giv \ff_s]
\en
and
\eq\label{p2-def}
   p_2(s,t) \Def \pr_{P_0}[\{\YYY_{P_0}(t) \cap \bC_{P}(t;2t)\ne \emptyset\}
         \cap \{\YYY_{P_0}(t) \cap \bC_{P'}(t;2t)\ne \emptyset\} \giv \ff_s],
\en
are close to probabilities
\eq
\begin{array}{rl}
   p_1^*(s,t) &:=\ \ex_{P_0}\{P^*[M\ui(t),M\ut(t)] \giv \ff_s\}; \\[1ex]
   p_2^*(s,t) &:=\ \ex_{P_0}\{P^*[M\ui(t),M\ut(t)] P^*[M\ui(t),M\uh(t)] \giv \ff_s\},
\end{array}\label{p-star-defs}
\en
that we
can approximate using Lemma~\ref{br-pr-intersect}.  Here, $M\ui$ is used to denote the
quantities~\Ref{M-defs} for~$\MBP_{P_0}$, coupled as above with~$\YYY_{P_0}$, and~$\ff_s$
to denote its history up to~$s$. $M\ut$ and~$M\uh$ are related to~$\bC_{P}(t;2t)$
and~$\bC_{P'}(t;2t)$ in similar fashion, through branching processes $\MBP_{P}$
and~$\MBP_{P'}$, which are independent of each other and of~$\MBP_{P_0}$.
Now the union~$\YYYs_{P_0}(t)$ of the islands
of~$\MBP_{P_0}(t)$ contains $\YYY_{P_0}(t)$, and the corresponding statement is true for
$\MBP_{P}(t)$ and~$\bC_{P}(t;2t)$ and for $\MBP_{P'}(t)$ and~$\bC_{P'}(t;2t)$.
From this, it follows that $p_1^*(s,t) \ge p_1(s,t)$ and that $p_2^*(s,t) \ge p_2(s,t)$.
We thus need only to show that the differences $p_1^*(s,t) - p_1(s,t)$ and
$p_2^*(s,t) - p_2(s,t)$ cannot be too large.
This we establish on the event~$\tA_{K,s} \in \ff_s$, defined as in~\Ref{ALs-def-cap} with
$M_i\ui(t)\l_0^i/i!$ for $H_i(t)$,  for a suitable choice of $K=K(\L)$.

The differences between $p_l^*(s,t)$ and~$p_l(s,t)$, $l=1,2$, arise from events on which an island
of a branching process~$\MBP$ makes an intersection that is not an intersection
in the corresponding $Y$-process.
Such an event can only occur if the island of~$\MBP$ is a ghost, or if the intersection
occurs at a part of an island that overlaps another island, so that an intersection may
have been counted twice using~$\MBP$.  Thus it will be enough to show that the probability
of there being a ghost or an overlapped island in~$\MBP_{P_0}(t)$ that
intersects~$\bC_{P}(t;2t)\cup\bC_{P'}(t;2t)$ is small, and that the same is true
for ghosts and overlapped islands of $\MBP_{P,P'}(t) := \MBP_{P}(t)\cup\MBP_{P'}(t)$
intersecting~$\YYY_{P_0}(t)$.

\begin{lem}\label{ghosts-etc}
Define $K(\L) := (40\log\L)^3$. Then,
on $\tA_{K(\L),s}$, for $t_{\L,x} := \frac1{2\l_0}\{\log\L + x\}$ as in~\Ref{t-lambda-def},
and with $|x| \le (1/6)\log\L$,
\[
  0 \Le  p_l^*(s,t) - p_l(s,t)
   \Eq O(e^x(1+e^x)\L^{-5/12}\{\log\L\}^{3d+12} + \L^{-1/6}\{\log\L\}^{d+6}),\quad l=1,2.
\]
\end{lem}

\proof
To make the necessary estimates, we begin by
coupling~$\MBP_{P_0}$ in $t\ge s$ to upper and lower processes
$\Xp$ and~$\Xm$ as in~\Ref{Xpm-bnds}, starting from the same state at time~$s$.
Define the event~$A'_+(K(\L),s)$ as in~\Ref{ALC-def}, but for the process~$\Xp$, and so
with $\l_0^+$ for~$\l_0$ and with $H^+$ for~$H$.
Then, on the event~$A'_+(K(\L),s) \cap A_{\th K,s}\ui$,
$M_0^+(t_{\L,x}) \le e\th K(\L)\L^{1/2} e^{x/2}$, and, in view of Theorem~\ref{H-at-t_L},
$A'_+(K(\L),s)$ has conditional probability at least $1 - c_2\L^{-8}$ on~$\tA_{K(\L),s}$.
Hence, 
the mean number of intersecting pairs of islands of~$\MBP_{P_0}$ at~$t_{\L,x}$ satisfies
\eq\label{mean-intersect}
   \ex \{N(t_{\L,x})I[A'_+(K(\L),s)] \giv \ff_s\cap \tA_{K(\L),s} \}
      \Le \half \{C  K(\L) \L^{1/2} e^{x/2}\}^2 p_\L^+ \Le C' e^x \{\log\L\}^{d+6},
\en
for some constants~$C,C'$, and also, from~\Ref{plp-def},
\eq\label{remainder-este}
   4\ex \{{M_0(t_{\L,x})}I[A'_+(K(\L),s)] \giv \ff_s\cap \tA_{K(\L),s}\}p_\L^+
        \Le C'' e^{x/2} \L^{-1/2} \{\log\L\}^{d+3}.
\en
Thus, taking $x = -(1/6)\log\L$, the conditional probability of any pair of islands
of~$\MBP_{P_0}$ intersecting
before time $(5/12)\l_0^{-1}\log\L$ is at most of order $O(\L^{-1/6}\{\log\L\}^{d+6})$
on $\tA_{K(\L),s}$, from~\Ref{mean-intersect}.  Also, if $x \le \half\log\L$,
the number of such intersecting pairs before time~$t_{\L,x}$ exceeds
$\log\L\{1 \vee e^x \{\log\L\}^{d+6}\}$
with conditional probability of order $O(\L^{-1/4}\{\log\L\}^{d+3})$,
in view of \Ref{mean-intersect}, \Ref{remainder-este}, Lemma~\ref{Poisson},
and the Chernoff inequalities for the Poisson distribution.

Thus, on~$\tA_{K(\L),s}$, except on an event of probability of order $O(\L^{-1/6}\{\log\L\}^{d+6})$,
all intersections of pairs of islands of~$\MBP_{P_0}$ occur after time~$(5/12)\l_0^{-1}\log\L$,
and, once more in view
of Theorem~\ref{H-at-t_L} applied to the dominating branching process~$\Xp$, they each give
rise to at most $C(\log\L)^{((d+1)\vee3)} \L^{1/12}e^{x/2}$
ghosts, except on an event of probability at most $c_a\L^{-7}$, since each such intersection
at worst gives rise to a ghost branching process starting with a single island of radius
$(3/2)\l_0^{-1}\log\L$.  Thus, on~$\tA_{K(\L),s}$, except on an event of conditional
probability of order
$O(\L^{-1/6}\{\log\L\}^{d+6})$, there are at most of order $O(\L^{1/12}\{\log\L\}^{2d+9}
e^{x/2}(1+e^x))$ islands in~$\MBP_{P_0}(t)$ whose intersections
with~$\bC_{P}(t;2t)\cup\bC_{P'}(t;2t)$ should be
discounted, and their radii cannot exceed $(3/2)\l_0^{-1}\log\L$.  Furthermore,
the number of islands in~$\MBP_{P}(t)$ and~$\MBP_{P'}(t)$ together is at most
$$
   C_aK(\L)e^{\l_0^+t(1+\e_{K(\L)})} \Eq O(\L^{1/2}\{\log\L\}^{3} e^{x/2}),
$$
except on an event of probability at most $c_a\L^{-8}$, by Theorem~\ref{H-at-t_L}(1).
Thus, on~$\tA_{K(\L),s}$, the
mean number of intersections that should be neglected is, off the exceptional events, at most of
order
\[
   O\Bl \L^{1/12}\{\log\L\}^{2d+9} e^{x/2}(1+e^x)\, \L^{1/2}\{\log\L\}^{3} e^{x/2} p_\L^+ \Br
     \Eq O(e^x(1+e^x)\L^{-5/12}\{\log\L\}^{3d+12}).
\]
This expectation, together with the complementary expectation from overlapping
islands and ghosts in $\MBP_{P,P'}$, bounds that part of the differences
$p_l^*(s,t)-p_l(s,t)$, $l=1,2$, arising off the exceptional events, and the exceptional events
together have probability of order~$O(\L^{-1/6}\{\log\L\}^{d+6})$.
The argument for this second expectation is the same, except that there is no conditioning,
making it equivalent to the previous argument with $s=0$, and with $A(K(\L),0)$ being
automatically satisfied because the inital state consists of just two singletons.
\ep

\begin{lem}\label{change-to-s}
For $M\ui,M\ut$ and~$M\uh$ as in Lemma~\ref{ghosts-etc}, and for $x \le \half\log\L$, we have
\eqa
   &&\left| \exs\biggl\{  \exp\Bigl\{-\sjd {d\choose j}  M_j\ui(t_{\L,x}) M_{d-j}\ut(t_{\L,x}) \l_0^d \L^{-1}
                 \Bigr\} \biggr\}  \right.\non \\
   &&\mbox{}\hskip2in\  \mbox{} \left.
      - \f\Bigl( (d+1)^{-2}\sjd {d\choose j} j!(d-j)! W_j\ui(s)  e^x \Bigr)  \right| \non\\[1ex]
   &&\qquad \ \Le Ce^x \{e^{-\b_{d+1}\l_0 s} (\log\L)^3  + (\log\L)^{4+\g_g}\L^{-\g_g/d} \}, \label{t-to-s}
\ena
and
\eqa
   &&\left| \exs\biggl\{  \exp\Bigl\{-\sjd {d\choose j}  M_j\ui(t_{\L,x})
            [M_{d-j}\ut(t_{\L,x})+M_{d-j}\uh(t_{\L,x})] \l_0^d \L^{-1}
                 \Bigr\} \biggr\}  \right.\non \\
   &&\mbox{}\hskip2in\  \mbox{} \left.
      - \biggl\{\f\Bigl( (d+1)^{-2}\sjd {d\choose j} j!(d-j)! W_j\ui(s)  e^x \Bigr)\biggr\}^2  \right| \non\\[1ex]
   &&\qquad \ \Le Ce^x (\log\L)^6\{e^{-\b_{d+1}\l_0 s}   + (\log\L)^{1+\g_g}\L^{-\g_g/d} \}, \label{t-to-s-2}
\ena
for any $s \le t_{\L,x}$,
where $\exs$ denotes expectation conditional on $\ff_s\cap \tA_{K(\L),s}$,
$W_j\ui(s) := (d+1)e^{-\l_0s}M_j\ui(s)\l_0^j/j!$, $\f$ is as in~\Ref{phi-def}, and~$C,C'$
are suitable constants.
\end{lem}

\proof
We give the proof of~\Ref{t-to-s-2}; that of~\Ref{t-to-s} is simpler.
For $t>s$, we bound $M^{(l)}$ above and below by $M^{(l)+}$ and~$M^{(l)-}$, $1\le l\le3$, as in the previous
lemma.  We then observe that, for $l=2,3$,
\eqs
  \lefteqn{M\ui_j(t_{\L,x})M^{(l)}_{d-j}(t_{\L,x})\l_0^d\L^{-1} }\\
       &\le& M\uip_j(t_{\L,x})M\ulp_{d-j}(t_{\L,x})(\l_0^+)^d\L^{-1} \\
       &=& (d+1)^{-2}j!(d-j)!e^{2\l_0^+t}W\uip_j(t_{\L,x})W\ulp_{d-j}(t_{\L,x}) \L^{-1} \\
      &\le& (d+1)^{-2}j!(d-j)!e^{(1+\h_\L)(\log\L + x)}W\uip_j(t_{\L,x})W\ulp_{d-j}(t_{\L,x}) \L^{-1} \\
      &\le& (d+1)^{-2}j!(d-j)!e^xW\uip_j(t_{\L,x})W\ulp_{d-j}(t_{\L,x}) \{1 + 2e\h_\L\log\L\},
\ens
where $\h_\L$ is as in~\Ref{eta-L-def}.
We now use Theorem~\ref{vector-control}\,(3) to deduce that
\eqa
\begin{array}{rl}
   \exs\Bl |W\uip_j(t_{\L,x}) - W\uip_j(s)| I[A'_+(K(\L),s)]  \Br
      &=\ O(K(\L)e^{-\b_{d+1}\l_0^+s})\phantom{XX} \\[2ex]
   \exs\Bl |W\ulp_{d-j}(t_{\L,x}) - W\ulp_{d-j}(\infty)| I[A'_+(K(\L),s)]  \Br
      &=\ O(K(\L)e^{-\b_{d+1}\l_0^+s}),   \phantom{XX}
\end{array} \label{W-diffs}
\ena
$l=2,3$, and Theorem~\ref{H-at-t_L}\,(2) to show that
$\pr[A'_+(K(\L),s) \giv  \ff_s \cap \tA_{K(\L),s}] \ge 1 - c_2\L^{-8}$.  It thus follows, also
using~ \Ref{expectation-bnd}, that
\eqa
   \lefteqn{\exs\biggl\{  \exp\Bigl\{-\sjd {d\choose j}  M_j\ui(t_{\L,x})
              [M_{d-j}\ut(t_{\L,x}) + M_{d-j}\uh(t_{\L,x})] \l_0^d \L^{-1}
                 \Bigr\} \biggr\}}\non \\
   &\ge& \exs\biggl\{  \exp\Bigl\{-(d+1)^{-2} d! \sjd    W_j\uip(s)
                     [W_*\utp(\infty) + W_*\uhp(\infty)] e^x \Bigr\} \biggr\}\non \\
   &&\  \mbox{}- O\bigl(e^x K(\L)^2\{e^{-\b_{d+1}\l_0^+s} + \h_\L\log\L\}  + \L^{-8}\bigr) \non\\
    &=& \Bigl\{\f\Bigl( (d+1)^{-2} d! \sjd   W_j\ui(s)  e^x \Bigr) \Bigr\}^2
              \label{lower-bnd-formula}\\
    &&\ \qquad \mbox{}  - O\bigl(e^x(\log\L)^6[e^{-\b_{d+1}\l_0s} + (\log\L)^{1+\g_g}\L^{-\g_g/d}]\bigr),  \non
\ena
since $W_j\uip(s) = W_j\ui(s)$ is $\ff_s$-measurable, and since $W_*\utp(\infty)$ and~$W_*\utp(\infty)$
are independent given~$\ff_s$, and each has the distribution of~$W_*(\infty)$
as in~\Ref{phi-def}.  The upper bound is proved in analogous fashion.
\ep

\begin{thm}\label{main-theorem}
   Choose $s_\L := (\a/2)\l_0^{-1}\log\L$ for some $0 < \a < 1/2$.  Let~$\b_r$
be as defined in Theorem~\ref{vector-control}, let $\tA_{K,s}$ be as for Lemma~\ref{ghosts-etc},
and define~$t_{\L,x}$ as in~\Ref{t-lambda-def}.
Then, for any $d$-dimensional gossip process satisfying~\Ref{volumes},
and for any $\g_1 < \g_0 := \min\{\third\a\b_{d+1},\twothirds\g_g/d\}$, there exists a
constant~$k_{\g_1} < \infty$ such that
\[
    \var\{L^{-1}V_{P_0}(2t_{\L,x}) \giv \ff_{s_\L} \cap \tA_{K(\L),s_\L}\} \Le k_{\g_1}\L^{-\g_1},
\]
uniformly in $|x| \le c_v\log\L$, where $c_v := \third\{\half \a\b_{d+1}\wedge \g_g/d\}$.
Furthermore, $\pr[\tA_{K(\L),s_\L}] \ge 1 - c_A\L^{-\g_2}$ for some $\g_2 > 0$ and~$c_A < \infty$.
\end{thm}

\proof
Starting from~\Ref{main-asymp-goss}, collecting the results of Lemmas \ref{br-pr-intersect},
\ref{ghosts-etc} and~\ref{change-to-s}, and taking $|x| \le \eighth\log\L$, we find
after comparing the various errors that
\eqa
   \lefteqn{\ex\{L^{-1}V_{P_0}(2t_{\L,x}) \giv \ff_{s_\L} \cap \tA_{K(\L),{s_\L}}\}
                 \Eq G(W\ui({s_\L}),x)}\non\\
     &&\mbox{}\hskip0.5in +
       O\bigl(e^x\{(\log\L)^3 \L^{-\a\b_{d+1}/2} + (\log\L)^{4+\g_g}\L^{-\g_g/d}\} +
                (\log\L)^{d+6+\g_g} \L^{-\g_d/d}\bigr),\non \\
    &&\ \Eq G(W\ui({s_\L}),x) + O\bigl(e^x(\log\L)^{4+\g_g} \L^{-3c_v}  +
                (\log\L)^{d+6+\g_g} \L^{-\g_g/d}\bigr),\label{vol-mean}
\ena
where
\eq\label{G-def}
    G(w,x) \Eq 1 - \f_d\Bigl(e^x(d+1)^{-2} d! \sjd w_j \Bigr).
\en
Similarly,
\eqa
  \lefteqn{\ex\{[L^{-1}V_{P_0}(2t_{\L,x})]^2 \giv \ff_{s_\L} \cap \tA_{K(\L),{s_\L}}\}}\non\\
       &&\Eq  \{G(W\ui({s_\L}),x)\}^2
           +  O\bigl(e^x(\log\L)^{4+\g_g} \L^{-3c_v}  +
                (\log\L)^{d+6+\g_g} \L^{-\g_g/d}\bigr).\phantom{XX}\label{vol-var}
\ena
Combining \Ref{vol-mean} and~\Ref{vol-var}, and bounding the error term, the
first statement of the theorem follows.

To bound $\pr[\tA_{K(\L),{s_\L}}]$ from below, we start with Theorem~\ref{H-at-t_L}\,(1).
Putting $s=0$ and then $t={s_\L}$, and setting $K=K(\L) = (40\log\L)^3$, it follows that
\[
   \pr[e^{-\l_0 {s_\L}}\|H({s_\L})\| \le C_a K(\L) e^{1/80}] \ \ge\ 1 - c_a\L^{-8},
\]
so that $\pr[A\ui_{\th K(\L),{s_\L}}] \ge 1 - c_a\L^{-8}$ if $\th := C_a e^{1/80}$. Then,
from Lemma~\ref{lower-bnd},
\[
   \pr[A\ut_{K(\L),{s_\L}}] \Eq \pr[H_{r(d)}({s_\L}) \ge K(\L)] \Eq \pr[\t^{r(d)}_{K(\L)} \le {s_\L}],
\]
bounded below by $1 - 2K(\L)\exp\{-{s_\L}\l_0/c_c\} = 1 - C\L^{-\g}$ for some $\g>0$,
because of the choice of~${s_\L}$.  Finally, from Lemma~\ref{PP-bound}\,(4),
since $K(\L)^{\shalf(1-\e_{K(\L)})}(2^{\e_{K(\L)}}-1) \asymp (\log\L)^{1/2}$ exceeds $42\log 2$
for all~$\L$ sufficiently large,
\eqs
   1 - \pr[A\uh_{K(\L),\e_{K(\L)},s_\L}] &\le& c_4
     \exp\{-{K(\L)}^{\shalf(1-\e_{K(\L)})}/28\} \\
   &\le& C_1 \exp\{-C_2(\log\L)^{3/2}\},
\ens
completing the proof.
\ep

\bigskip
This theorem is the basis for the main result of the paper, showing that the distribution of
the path~$L^{-1}V_{P_0}(t)$ is concentrated close to its conditional mean.  To complete
the proof, we need first to have an \adb{expression for the conditional mean.}

\begin{lem}\label{mean-lemma}
Uniformly in $|x| \le c_v\log\L$, and with
$s_\L := (\a/2)\l_0^{-1}\log\L$ for some $0 < \a < 1/2$,
\[
   \ex\Blb L^{-1}V_{P_0}(\l_0^{-1}\{\log\L + x\})
        \giv \ff_{s_\L} \cap \tA_{K(\L),{s_\L}} \Brb \Eq h_d(x  + \log C_d + \log W_*\ui(s_\L)) + O(\L^{-\g_1}),
\]
where~$h_d$ is as in~\Ref{h-eqn}, $C_d := (d+1)^{-1}d!$,
and $\g_1$ is as for Theorem~\ref{main-theorem}. Also, if~$P$ is independently and uniformly
chosen on~$C$, the time~$\t_P := \inf\{t\ge0\colon\, P \in \YYY_{P_0}(t)\}$ satisfies
\[
   \pr[\t_P > \l_0^{-1}(\log\L + x)] \Eq
    \ex\Bigl(\exp\Bigl\{-e^x C_d W_*\ui(\infty) W_*\ut(\infty)
      \Bigr\}\Bigr)  + O(\L^{-\g_3}),
\]
for some~$\g_3 >0$, where $W_*\ui(\infty)$ and~$W_*\ut(\infty)$ are independent.
\end{lem}

\proof
By arguing as for Lemma~\ref{change-to-s}, using Theorem~\ref{vector-control}(3),
it follows that, uniformly for~$x$ as in Theorem~\ref{main-theorem},
\eqa\label{vol-mean-new}
  \lefteqn{ \ex\{L^{-1}V_{P_0}(2t_{\L,x}) \giv \ff_{s_\L} \cap \tA_{K(\L),{s_\L}}\}}\non\\
     &=& \ex\{G(W\ui({s_\L}),x)\giv \ff_{s_\L} \cap \tA_{K(\L),{s_\L}}\}
                     + O(\L^{-\g_1}) \non\\
     &=&   \ex\{G(W_*\ui({s_\L})\bone,x)\giv \ff_{s_\L} \cap \tA_{K(\L),{s_\L}}\}
            + O(\L^{-\g_1}) \non\\
     &=& 1 - \f_d(C_d e^{x} W_*\ui({s_\L})) + O(\L^{-\g_1}),
\ena
where~$C_d$ is as defined above, and the first part follows from the definition
of~$h_d$ following~\Ref{phi-def}.   The final result follows from taking the unconditional
expectation in~\Ref{vol-mean-new}, and applying Theorem~\ref{vector-control}\,(3):
\eqs
   1 -  \ex\{L^{-1}V_{P_0}(2t_{\L,x})\} &=& \ex\{\f_d(C_de^{x} W_*\ui(s_\L))\}
      + O(\L^{-\g_1}) + O(\L^{-\g_2}) \\
    &=& \ex\{\f_d(C_de^{x} W_*\ui(\infty))\} + O(\L^{-\g_3}),
\ens
for some $\g_3 >0$.
\ep

Note that the form of the neighbourhoods~$\KK$
only comes into the formulae of Lemma~\ref{mean-lemma} through their volume~$v(\KK)$,
which is implicitly present in the time scaling by~$\l_0$ in the definition of~$t_{\L,x}$.

We are now in a position to prove the pathwise approximation to~$L^{-1}V_{P_0}(t)$.
Before doing so, we note that $h_d$ is the distribution function of a sum of independent
random variables $Z_1$ and~$Z_2$, where $-Z_1$ has a standard Gumbel distribution,
and~$Z_2$ is distributed as $-\log W_*(\infty)$.  This is because $h_d(x) := 1-\f_d(e^{x})$ can
be rewritten in alternative form as
$$
    1-\ex\Blb e^{-e^Z} \Brb  \Eq  \pr[-Z_1 \ge -Z] \Eq \pr[Z_1+Z_2 \le x],
$$
where $Z := x  + \log W_*(\infty)$ is independent of~$Z_1$.

\begin{thm}\label{path-theorem}
For any $d$-dimensional gossip process satisfying \Ref{volumes},
there exists a random variable~$U$ such that
\[
   \pr[\sup_x|L^{-1}V_{P_0}(\l_0^{-1}\{\log\L + x\}) - h_d(x + \log C_d + U)| > 4\L^{-a_1}] \Le c\L^{-a_2},
\]
for some positive $a_1,a_2$ and $c<\infty$.  Here, $e^{U}$ has the distribution
of~$W_*^{[d+1]}(\infty)$ as in Remark~\ref{W-rk},
whose Laplace transform satisfies~\Ref{phi-eqn}, $C_d$ is as defined in Lemma~\ref{mean-lemma},
$h_d$ is as in~\Ref{h-eqn}, and~$\L$ is as defined in~\Ref{Lambda-def}.
\end{thm}

\proof
Take $s_\L := (\a/2)\l_0^{-1}\log\L$ for some $0 < \a < 1/2$.
It then follows from Theorem~\ref{main-theorem}, Lemma~\ref{mean-lemma} and
Chebyshev's inequality that, for any $a < \g_1/2$,
\eqa\label{point-appxn}
  \lefteqn{\pr[|L^{-1}V_{P_0}(\l_0^{-1}\{\log\L + x\}) - h_d(x  + \log C_d + \log W_*\ui(s_\L))| > 2\L^{-a}
       \giv \ff_{s_\L} \cap \tA_{K(\L),{s_\L}} ]}\non\\
  &&\qquad \Le k_{\g_1}\L^{2a-\g_1} \phantom{XXXXXXXXXXXXXXXXXXXXXXXXXXXXXXX}
\ena
uniformly for $|x| \le c_v\log\L$.

Now, representing~$h$ using $Z_1$ and~$Z_2$ as above, we have, for $x>0$,
\eqs
     h_d(-x) &\le& \pr[-Z_1 \ge x/2] + \pr[\log W_*(\infty) \ge x/2] \Le 2e^{-x/2};\\
     1- h_d(x) &\le& \pr[-Z_1 \le -x/2] + \pr[\log W_*(\infty) \le -x/2] \Le \exp\{-e^{x/2}\} + \e(x),
\ens
where $\e(x)$ goes to zero super-exponentially fast, by \Ref{W-lower-tail-bnd}.
Thus, for $x_\L^+ := \half c_v\log\L$,
$1-h_d(x_\L^+) = O(\L^{-1})$, and, for $x_\L^- := -\half c_v\log\L$, $h_d(x_\L^-) \Le 2\L^{-c_v/4}$. Now
choose an increasing sequence of~$m_\L := \lceil \L^b \rceil$ points $x_\L^{(j)}$ between $x_\L^-$ and~$x_\L^+$,
with $x_\L^{(0)} = x_\L^-$ and $x_\L^{(m_\L)} = x_\L^+$, in such a way that
$h_d(x_\L^{(j+1)}) - h_d(x_\L^{(j)}) \le \L^{-b}$ for each~$j$.  On~$\tA_{K(\L),s_\L}$,
from~\Ref{ALs-def},
\eq\label{log-W-tails}
   -\log W_*\ui(s_\L) \Le \l_0 s_\L \Eq \half\a\log\L;\qquad \log W_*\ui(s_\L) \Le \log (dK(\L)),
\en
so that, choosing $\a \le \half c_v$, $|\log W_*\ui(s_\L)| \le \quarter c_v\log\L$ for all $\L$ large enough.
Thus we can take $x = x_\L^{(j)} - \log W_*\ui(s_\L)$ in~\Ref{point-appxn} for each $0\le j\le m_\L$,
obtaining
\eqa
  \lefteqn{\pr\Bigl[\max_{0\le j\le m_\L}|L^{-1}V_{P_0}(\l_0^{-1}\{\log\L + x_\L^{(j)} - \log W_*\ui(s_\L)
       - \log C_d\} - h_d(x_\L^{(j)})|}  \label{many-point-appxn}\\
    &&\qquad\mbox{}\hskip2.2in  > 2\L^{-a} \Giv \ff_{s_\L} \cap \tA_{K(\L),{s_\L}} \Bigr]
   \Le 2k_{\g_1}\L^{2a+b-\g_1}, \non 
\ena
from which it follows directly, because both $V_{P_0}(t)$ and~$h_d(x)$ are non-decreasing in their arguments,
that
\eqa
 \lefteqn{\pr\Bigl[\sup_{x}|L^{-1}V_{P_0}(\l_0^{-1}\{\log\L + x \}) -h_d(x + \log C_d + \log W_*\ui(s_\L))|}
     \label{all-points-appxn-*} \\
   &&\qquad\mbox{}\hskip2in > 2\L^{-a} + \L^{-b}  \Giv \ff_{s_\L} \cap \tA_{K(\L),{s_\L}} \Bigr]
   \Le k\L^{2a+b-\g_1}, \non
\ena
for a suitable constant~$k$.

The above estimate almost completes the theorem; it remains to replace $W_*\ui(s_\L)$
by $W_*\uip(\infty)$, where $W\uip$ is as in Lemma~\ref{ghosts-etc}, and then to remove the
conditioning.  For the former, we have
\[
    \ex\{|W_*\ui(s_\L) - W_*\uip(\infty)|I[A'_+(K(\L),s_\L)] \giv \ff_{s_\L} \cap \tA_{K(\L),{s_\L}}\}
      \Eq O(K(\L)\exp\{-\b_{d+1}\l_0^+ s_\L\}),
\]
from Theorem~\ref{vector-control}(3), together with
$\pr[A'_+(K(\L),s_\L) \giv \ff_{s_\L} \cap \tA_{K(\L),{s_\L}}] \ge 1-c_2\L^{-8}$,
from Theorem~\ref{H-at-t_L}(2). Thus
\[
   \pr[|W_*\ui(s_\L) - W_*\uip(\infty)| > \L^{-b} \giv \ff_{s_\L} \cap \tA_{K(\L),{s_\L}}]
     \Eq O(\L^{-b'}),
\]
for any $b' < \half\b_{d+1}\a - b$, and thus for any $b' < 3c_v  - b$, and hence also
\[
   \pr[|\log W_*\ui(s_\L) - \log W_*\uip(\infty)| > \L^{-b/2} \giv \ff_{s_\L} \cap \tA_{K(\L),{s_\L}} \cap
        \{W_*\uip(\infty) \ge \L^{-b/2} \} ]  \Eq O(\L^{-b'}).
\]
Since~$h_d$ has density bounded by~$1/e$, the maximum of the density of the Gumbel distribution,
it follows that
\eqa
  \lefteqn{\pr\Bigl[\sup_{x}|L^{-1}V_{P_0}(\l_0^{-1}\{\log\L + x \})
                  -h_d(x + \log C_d + \log W_*\uip(\infty))|} \non\\
      && \qquad\qquad > (2\L^{-a} + \L^{-b} + e^{-1}\L^{-b/2})
     \Giv\ff_{s_\L} \cap \tA_{K(\L),{s_\L}} \cap \{W_*\uip(\infty) \ge \L^{-b/2}\}\Bigr] \non\\
        && \Le k\L^{2a+b-\g_1} + k'\L^{-b'}.\phantom{XXX} \label{all-points-appxn-c}
\ena
But now, from~\Ref{W-lower-tail-bnd},
$\pr[W_*\uip(\infty) < \L^{-b/2}] = o(\L^{-u})$ for any $u>0$, and $\pr[\tA_{K(\L),{s_\L}}] \ge
1 - c_A\L^{-\g_2}$ from Theorem~\ref{main-theorem}, so that
\eqa
  \lefteqn{\pr\Bigl[\sup_{x}|L^{-1}V_{P_0}(\l_0^{-1}\{\log\L + x \})
                  -h_d(x + \log C_d + \log W_*\uip(\infty))| > 2\L^{-a} + \L^{-b} + e^{-1}\L^{-b/2} \Bigr]} \non\\
       && \qquad  \Le k\L^{2a+b-\g_1} + c_A\L^{-\g_2} + c'\L^{-b'}.
      \phantom{XXXXXXXXXXXXXXXXXXXX}\label{all-points-appxn}
\ena
Now take $2a=b=(c_v \wedge \{\g_1/4\})$ and $b' = c_v$ to complete the proof,
with $a_1 = b/2$ and $a_2 = \min\{\g_1/2,\g_2,c_v\}$.
\ep

For small worlds processes, the counterpart of Theorem~\ref{main-theorem} can be proved
in entirely similar fashion.  Lemma~\ref{mean-lemma} is also correct, if $C_d$
is replaced by
\eq\label{tC-def}
   \tC_d \Def d^{-2}d!\{(d+1)-1\} \Eq (d-1)!\,;
\en
the change  reflects both the difference in~$r(d)$ between
gossip and small world processes, and the subtracted component arising from intersections forbidden
in the small world context: see \cite{BR}, Section~4 for more details.  Note that the
expression in \cite{BR}, Theorem~4.2, appears different from that obtained here.  It is,
however, the same, with the transformations $r\to d$, $2x \to x$, $\a_l \to v(\KK){d\choose l}$,
$W \to d^{-1}W_*(\infty)$, and noting that $L\r = \L / d!$.  For the limiting random
variable~$W_*(\infty)$, the Poisson process of descendants of the first individual
now has intensity
\[
   \r v(\KK) d u^{d-1} \Eq \l_0^d u^{d-1}/(d-1)!,\qquad u > 0,
\]
giving
\eq\label{phi-eqn-sw}
   \f(\th)
     \Eq  \exp\Blb - \int_0^\infty \frac1{(d-1)!}\{1 - \f(\th e^{-x})\}x^{d-1}\,dx \Brb,
\en
the same equation as~\Ref{phi-eqn} for gossip processes, except that now $d$ is replaced by~$d-1$;
in view of~\Ref{time-scaled}, this is not surprising.  We thus have the corresponding
law of large numbers approximation.


\begin{thm}\label{path-theorem-2}
For any $d$-dimensional small world process satisfying \Ref{volumes},
there exists a random variable~$\tU$ such that
\[
   \pr[\sup_x|L^{-1}V_{P_0}(\l_0^{-1}\{\log\L + x\}) - h_{d-1}(x + \log \tC_d + \tU)| > 4\L^{-\ta_1}]
      \Le \tc\L^{-\ta_2},
\]
for some positive $\ta_1,\ta_2$ and $\tc < \infty$.   Here,  $e^{\tU}$ has the distribution of~$W_*^{[d]}(\infty)$
as defined in Remark~\ref{W-rk}, $\tC_d$ is as defined in~\Ref{tC-def},
$h_r$ is as in~\Ref{h-eqn}, and~$\L$ is as defined in~\Ref{Lambda-def}.
\end{thm}

For the gossip process studied by Chatterjee \& Durrett~\cite{ChatDurr}, $C$ is an $N\times N$ torus, so that
$d=2$ and $L=N^2$, and $\KK(P,s) = B(P,s/\sqrt{2\p})$ is a Euclidean ball, so that $v(\KK) = 1/2$
and $c_g=0$.  They take $\r = N^{-\a}$ for any $\a<3$, so that $\l_0 = (d!\r v(\KK))^{1/(d+1)}
= N^{-\a/3}$ and $\L = L\l_0^d/v(\KK) = 2N^{2(1-\a/3)}$; also $C_d = 2/3$ in Lemma~\ref{mean-lemma}.
Then the pathwise approximation
given in Theorem~\ref{path-theorem} matches $N^{-2}V(N^{\a/3}[2(1-\a/3)\log N + x + \log2])$
with $h(x+\log(2/3) + U)$.  This is seen to be the same as in Chatterjee \& Durrett~\cite{ChatDurr},
noting that $U$ has the same distribution as their~$\log M$, and that our choice of~$h(\cdot)$ as
a solution of~\Ref{h-eqn} corresponds to their $h(\cdot + \log3)$, since they implicitly
choose the solution of~\Ref{h-eqn} that has $\lim_{t\to-\infty}e^{-t}h'(t) = 1/3$.

\subsection{Complete coverage}\label{coverage}
Theorem~\ref{path-theorem} and~\Ref{vol-mean-new}, together with the fact
that $2t_{\L,x} = \l_0^{-1}\{\log\L + x\}$, show that only a negligible fraction
of~$C$ is covered at times much before $\l_0^{-1}\log\L$, with~$C$ then becoming
essentially covered in deterministic fashion on a time scale of order~$O(1)$ around
$\l_0^{-1}\log\L$.  Randomness is only visible in the time shift of $-\l_0^{-1}U$
in the origin of the transition from uncovered to covered.

One can also, as in Chatterjee \& Durrett~\cite{ChatDurr}, consider how long it takes
until~$C$ is entirely covered.  We show that complete coverage is achieved
relatively soon after time $\l_0^{-1}\log\L$, under the assumption that, for each~$s$,
$C$ can be covered by~$n(s)$ islands of the form~$\KK(P,s)$, where~$n(s)$ satisfies
\eq\label{coverage-assn}
    n(s) \Le c_0 L/\{v(\KK)s^d\},\qquad 0 < s < L^{1/d},
\en
for some $c_0$.

\begin{thm}\label{coverage-thm}
If~\Ref{coverage-assn} is satisfied for a $d$-dimensional gossip process
satisfying~\Ref{volumes}, then, except on a set of probability of order
$O(\L^{-\d})$, for some $\d>0$,  the whole of~$C$ is covered before time
$$
   \t(\L,s) + \frac2{\l_0}\Blb \frac{72 \log\L }{d!} \Brb^{1/d},
$$
where $\t(\L,s) = \l_0^{-1}\{\log\L + O((\log\L)^{1/(d+1)})\}$. The
corresponding result holds for a small world process, now with
$\t(\L,s) = \l_0^{-1}\{\log\L + O((\log\L)^{1/d})\}$.
\end{thm}

\proof
In a gossip process, in view of Theorem~\ref{path-theorem}, at least
$\quarter \L/d!$ contacts have been made up to
time~$\t(\L,s) := \l_0^{-1}\{\log\L + \mm - U - \log C_d\} + \l_0^{-1}$,
with probability $1 - O(\L^{-\d})$ for some $\d>0$,
where~$\mm$ denotes the median $h_{d}^{-1}(1/2)$; this is because, at this time,
a volume of about~$L/2$ has been generating contacts for a time interval of
at least~$\l_0^{-1}$.
Write $\ps := d!/4$, and note also that, from~\Ref{W-lower-tail-bnd}, $- U \le (k\log\L)^{1/(d+1)}$
except on an event of probability of order $O(\L^{-8})$, if~$k$ is chosen large enough.

Cover~$C$ with islands of radius $s =  k\l_0^{-1}$; with~$c_0$ as in~\Ref{coverage-assn}  this can be done using
at most $ c_0 L/\{v(\KK)(k/\l_0)^d\}$ islands.
Then, recalling \Ref{volumes},  the
probability that any of these islands contains none of $P_1,\ldots, P_{\lceil \ps\L \rceil}$
is at most
\[
    \frac{ c_0v(\KK)\L }{v(\KK)k^d}\,\Bl 1 - \frac{v(\KK)k^d}{2v(\KK)\L }\Br^{\lceil \ps\L \rceil}
      \Le \frac{c_0\L }{k^d}\,e^{-\ps k^d/2}\ =:\ \p_k,
\]
if~$\L$ is large enough that $c_g(k \L^{-1/d})^{\g_g} \le 1/2$.
On the complementary event, all of the islands, and thus all of~$C$, are covered after
an additional time of at most~$2k\l_0^{-1}$.
So take $k = \{(18/\ps)\log\L \}^{1/d}$, to make
$$
    \p_k \Le \frac{c_0\ps}{18\log\L }\L^{-8}.
$$

In a small world process, we apply Theorem~\ref{path-theorem-2} instead.
With probability $1 - O(\L^{-\d})$ for some $\d>0$, at least~$\quarter L\r = \quarter \L/d!$
shortcuts have been encountered by time~$\t(\L,s) := \l_0^{-1}\{\log\L + \tilde\mm - \tU -
\log \tC_d\}$, where~$\tilde\mm$ denotes the median $h_{d-1}^{-1}(1/2)$,
since then about half of~$C$ has been covered. Because~$e^{\tU}$ has the
distribution of~$W_*^{[d]}(\infty)$ rather than that of~$W_*^{[d+1]}(\infty)$,
we have $- \tU \le (k\log\L)^{1/d}$, except on an event of probability of
order~$O(\L^{-8})$, if~$k$ is chosen large enough, in view of~\Ref{W-lower-tail-bnd}.
The remainder of the argument is as for the gossip process.
\ep

\subsection{Manifolds with boundary}\label{manifolds}
The assumption that the manifold~$C$ is homogeneous simplifies the argument substantially.
However, the adjustments needed if~$C$ is taken to be a `reasonable' finite subset of a
homogeneous manifold, such as a rectangle in~$\re^2$ or a spherical cap, are not great.
The principal requirement is that most islands do not intersect the boundary~$\partial C$.
Let $C_\d := \{P\in C\colon\, \KK(P,\d) \cap \partial C \ne \emptyset\}$ denote the
$\d$-neighbourhood of the boundary of~$C$, and assume that its volume is not too large
when compared with that of~$C$:
\eq\label{annulus-assn}
    L^{-1}|C_\d| \Le c_b \d (v(\KK)/L)^{1/d}, \qquad 0 < \d \le 2\log\L,
\en
for some constant~$c_b$. For instance, with~$C$ a square of side~$\sqrt L$ in~$\re^2$
and with $\KK(P,s)$ a disc of radius~$s$, $d=2$ and $v(\KK) = \pi$, and
$L^{-1}|C_\d| \Le 4\d L^{-1/2}$ satisfies~\Ref{annulus-assn} with $c_b = 4/\sqrt{\pi}$.
As is clear from the preceding argument, only times
less than $\d_\L := 2\l_0^{-1}\log\L$ play a significant part, and the probability of an island
with randomly chosen centre $P \in C$ intersecting~$\partial C$ before time $\d_\L$
is then at most
\eq\label{boundary-prob}
    c_b\d_\L\,(v(\KK)/L)^{1/d} \Eq 2c_b \L^{-1/d} \log\L
\en
under the assumption~\Ref{annulus-assn}.

In order to make the arguments of Section~\ref{calculations} work, it is enough to be able
to bound the growth process above and below by branching processes with constant growth
rates $\l_0^+$ and~$\l_0^-$, respectively, which are close enough to one another,
and to have~\Ref{intersection-prob} hold,
with an error estimate similar to that in~\Ref{intersection-error}.  For the latter,
if~\Ref{volumes} is satisfied, then~\Ref{boundary-prob} implies that
\eq\label{intersection-prob-bdy}
   \ppp_L(t,u)  \Eq L^{-1}v(\KK)(t+u)^d(1 + R_L(t,u)),
\en
where
\eq\label{intersection-error-bdy}
   |R_L(t,u)| \Le c_g\{(t + u)(v(\KK)/L)^{1/d}\}^{\g_g} + c_b\max\{t,u\}\,(v(\KK)/L)^{1/d},
\en
so that all that is needed is to replace the exponent $\g_g$ by $\tilde\g_g := \min\{\g_g,1\}$
and the constant~$c_g$ by $\tilde c_g := c_g+c_b$ in~\Ref{intersection-error},
when making intersection calculations.  For the former, the upper bound~$\l_0^+$ given
in~\Ref{lambda-0-pm} still holds.  A lower bound is obtained by neglecting
any contacts to points of~$C_{\d_\L}$, and taking
\eq\label{tilde-lambda}
   \tilde\l_0^- \Def \l_0^-(1 - L^{-1}|C_\d|) \ \ge\ \l_0^-\{1  - 2c_b\L^{-1/d}\log\L\},
\en
by~\Ref{annulus-assn} and~\Ref{boundary-prob}.  Thus, once again, the previous arguments can be
carried through, if, in the definition~\Ref{eta-L-def}, $\g_g$ is replaced by~$\tilde\g_g$, and
$c_g$ by $c_g + 4dc_b/3$.  This leads to the following result.

\begin{thm}\label{bdy}
For a $d$-dimensional gossip process on a finite subset~$C$ of a
homogeneous manifold that satisfies \Ref{volumes} and~\Ref{annulus-assn}, the conclusion
of Theorem~\ref{path-theorem} holds, but with different constants $a_1,a_2$ and~$c$.
For small worlds processes, under the same assumptions,
the conclusion of Theorem~\ref{path-theorem-2} holds, again with different constants.
\end{thm}

For the result corresponding to Theorem~\ref{coverage-thm}, it is also necessary to make the
explicit assumption in connection with~\Ref{coverage-assn}, which was previously guaranteed for all~$\L$
such that $c_g(k \L^{-1/d})^{\g_g} \le 1/2$,
that each of the $n(s)$ sets~$\KK(P,s)$ used to cover~$C$ satisfies
\eq\label{big-overlap}
     |\KK(P,s) \cap C| \ \ge\ \half s^d v(\KK).
\en
In the current context, this is no longer automatic, because part of a set $\KK(P,s)$ may
lie outside~$C$.  The proof otherwise runs without any change, and the conclusion
of Theorem~\ref{coverage-thm} holds under this extra assumption. Note also that, under
the extra assumptions \Ref{annulus-assn} and~\Ref{big-overlap} of this section, the
manifold~$C$ could also be allowed to consist of a number of disconnected components.
In particular, the requirement \Ref{big-overlap} applied with $s = k\l_0^{-1}$ would
prevent Theorem~\ref{coverage-thm}
from being justified if there were components that were too small, and this is to be
expected, since it may take an extremely long time for a very small component of~$C$ to be
hit by a sequence of randomly chosen points of~$C$.

\section*{Acknowledgement}
ADB thanks the Institute for Mathematical Sciences at the National University of Singapore,
and the mathematics departments of the University of Melbourne and Monash University, for
their kind hospitality while part of the work was undertaken. GDR also thanks the
Institute for Mathematical Sciences at the National University of Singapore for their kind hospitality.
The authors are grateful to Dirk Schl\"{u}ter for helpful discussions, and to a referee,
whose detailed input has greatly improved the presentation of the paper.

\end{document}